\documentclass[final]{elsarticle}
\usepackage{amsmath, amsthm, amssymb}
\usepackage{amsfonts}
\usepackage{mathrsfs}
\usepackage{etoolbox}
\usepackage{graphicx}
\usepackage[unicode=true,bookmarks=true, colorlinks=true]{hyperref}
\hypersetup{pdfauthor={hyu}}
\newcommand{\cmt}[1]{{}}  % define a quick way to comment out text

\newcommand{\nn}{\nonumber}
\newcommand{\veps}{\varepsilon}
\DeclareGraphicsExtensions{.pdf,.pdf,.png,.jpg}

\newtheorem{theorem}{Theorem}
\newtheorem{lemma}{Lemma}
\newtheorem{remark}{Remark}
\usepackage[margin=1.2in]{geometry}  % use 1cm for draft

\begin{document}
\begin{frontmatter}

\title{Energetic Spectral-Element Time Marching Methods
    for Phase-Field Nonlinear Gradient Systems}

\author[label2,label3]{Shiqin Liu}
\ead{sqliu@lsec.ac.cc.cn}

\author[label3,label2]{Haijun Yu\corref{cor1}}
\ead{hyu@lsec.cc.ac.cn}

\address[label2]{School of Mathematical Sciences,
	University of Chinese Academy of Sciences, Beijing 100049, China}
\address[label3]{LSEC \& NCMIS, Institute of
	 	Computational Mathematics and Scientific/Engineering
	 	Computing, Academy of Mathematics and Systems Science,
	 	Beijing 100190, China}
\cortext[cor1]{Corresponding author}

\begin{abstract}
   We propose two efficient energetic spectral-element methods in time
   for marching nonlinear gradient systems with the phase-field Allen--Cahn
   equation as an example: one fully implicit nonlinear method and one
   semi-implicit linear method. Different from other spectral methods
   in time using spectral Petrov-Galerkin or weighted Galerkin
   approximations, the presented implicit method employs an energetic
   variational Galerkin form that can maintain the mass conservation and energy
   dissipation property of the continuous dynamical system. Another advantage
   of this method is its superconvergence. A high-order extrapolation is adopted for the
   nonlinear term to get the semi-implicit method. The semi-implicit
   method does not have superconvergence, but can be improved by a few
   Picard-like iterations to recover the superconvergence of the implicit
   method. Numerical experiments verify that the method using Legendre elements
   of degree three outperforms the 4th-order implicit-explicit backward
   differentiation formula and the 4th-order exponential time difference
   Runge-Kutta method, which were known to have best performances in solving
   phase-field equations. In addition to the standard Allen--Cahn equation, we
   also apply the method to a conservative Allen--Cahn equation, in which the
   conservation of discrete total mass is verified. The applications of the
   proposed methods are not limited to phase-field Allen--Cahn equations. They
   are suitable for solving general, large-scale nonlinear dynamical systems.
\end{abstract}

% REQUIRED

%% REQUIRED
%\begin{MSCcodes}
%    65M70, 65L05, 65M12
%\end{MSCcodes}
%
\end{frontmatter}

\section{Introduction}

The phase-field method is a highly effective computational tool for modeling
and predicting morphological evolutions in
material science~\cite{cahn_free_1958, allen_microscopic_1979, chen_phase_2002,
greenwood_free_2010}, fluid dynamics~\cite{anderson_diffuseinterface_1998,
lowengrub_quasiincompressible_1998, yue_diffuse_2004, qian_variational_2006},
biological systems~\cite{du_phase_2004, oden_general_2010,
wise_threedimensional_2008, miranville_long_2019}, etc. However, the small parameter
of interface thickness in phase-field models makes them numerically challenging
stiff and nonlinear problems. In this paper, we propose efficient energetic
spectral-element time marching methods for phase-field gradient systems, using
the following Allen--Cahn equation as an illustrative example:
\begin{align}
	\label{eq:AC}
   u_t - \varepsilon\Delta u + \frac1\varepsilon f(u) & = 0, \quad x
   \in \Omega\times[0,T], \\ u|_{\partial{\Omega}} &=
   0, \label{eq:DirichletBC}\\ u(x, t) &= u_0(x), \quad \text{at }
   t=0, \label{eq:IV}
\end{align}
where $f(u)=F'(u)$ and $F(u)$ represents a double-well potential, which is
typically chosen as
\begin{equation}
    \label{eq:double-wellF}
   F(u) = \frac{1}{4} (1-u^2)^2.
\end{equation}
For simplicity, we only analyze the Dirichlet condition \eqref{eq:DirichletBC}.
It is straightforward to extend the results in this paper to other types of
boundary conditions.

Equation \eqref{eq:AC} is introduced by Allen and
Cahn~\cite{allen_microscopic_1979} to describe the motion of
anti-phase boundary motion in crystalline solids. An important feature of the
Allen--Cahn equation is its energy dissipation property: The equation can be
viewed as a $L^2$ gradient flow of the Ginzburg-Landau free
energy~\cite{cahn_free_1958,liu_phase_2003}, which is defined as
\begin{equation}\label{eq:GL}
  E(u):= \int_{\Omega} \bigl(\frac{\varepsilon}{2}|\nabla u|^2 +
  \frac{1}{\varepsilon}F(u)\bigr)\, dx.
\end{equation}
By pairing (\ref{eq:AC}) with $u_t$ or $-\varepsilon\Delta u +
\frac{1}{\varepsilon} f(u)$, we find the energy law for (\ref{eq:AC}):
\begin{equation*}
    \frac{\partial}{\partial t}E(u(t)) =
    - \int_{\Omega} | u_t|^2 dx
    = - \int_{\Omega}
    |-\varepsilon\Delta u + \frac{1}{\varepsilon}f(u)|^2\ dx
    \le 0.
\end{equation*}
When solving phase field problems, it is important to verify the algorithms'
energy stability property at the discrete level, making energy-stable schemes
highly preferable. There are several popular numerical techniques for achieving
energy stability for large time step sizes, which can be roughly classified
into two categories. The first category maintains original energy (with a
possible addition of a diminishing momentum term) dissipation, including
methods such as the secant method~\cite{du_numerical_1991}, convex splitting
methods~\cite{elliott_global_1993,eyre_unconditionally_1998,wang_energy_2011,
shen_second_2012}, and linearly stabilization
schemes~\cite{chen_applications_1998,zhu_coarsening_1999, xu_stability_2006,
shen_numerical_2010, wang_efficient_2018, wang_energy_2020, li_double_2023},
ETD or IMEX Runge-Kutta methods\cite{guo_efficient_2016,wang_efficient_2016,
zhu_fast_2016,shin_unconditionally_2017, zhang_explicit_2021, ju_maximum_2021,
fu_energy_2022a}. The methods in the second category maintain stability of
modified energies with auxiliary variables. For example, in invariant energy
quadratization (IEQ)~\cite{yang_linear_2018, yang_efficient_2018a,
yang_convergence_2020}, the bulk potential is transformed into a quadratic form
using a set of new variables, and the nonlinear terms are semi-explicitly
treated, resulting in linear and unconditionally stable systems. This approach
is a generalization of augmented Lagrange multiplier(ALM)
method~\cite{guillen_linear_2013,guillen_second_2014}. The scalar auxiliary
variable (SAV) methods~\cite{shen_new_2019, shen_scalar_2018,
shen_convergence_2018} are built on a similar methodology, but are more
efficient by introducing an auxiliary scalar variable (instead of an auxiliary
function in IEQ approach). SAV methods can achieve second-order unconditionally
stability and only requires solving linear, decoupled systems with constant
coefficients at each time step.

However, designing efficient high order ($\ge 3$) schemes that can keep the
mass conservation in machine accuracy and keep energy dissipation with
reasonable time step sizes is still a challenging task. In this paper, we
introduce an implicit and a semi-implicit spectral-element time marching
schemes based on an energetic weak formulation. Let $H^1:=W^{1,2}(\Omega)$ with
$W^{k,p}(\Omega)$ denoting standard Sobolev space. Here $\Omega \subset
\mathbb{R}^d(d = 2, 3)$ is a bounded domain with $C^{1,1}$ boundary $\partial
\Omega$ or a convex polygonal domain. We define
\begin{equation*}
   V_{[s,t]}:= H^1\bigl(s,t; H_0^1(\Omega)\bigr), \quad
   V^{u_0}_{[s, t]}:= \bigl\{ u\in V_{[s, t]} \mid u(x,s) = u_0(x) \bigr\},
\end{equation*}
where $H_0^1(\Omega) = \{\, u \in H^1(\Omega),\ u|_{\partial \Omega} = 0\,\}$.
The energetic weak form is presented as follows:
\begin{equation}
   \label{eq:ACweakform}
   \mbox{Find}\ u\in V^{u_0}_{[0,T]},\ \text{s.t.}\ \int_0^T ( u_t, v_t
   ) + \varepsilon( \nabla u, \nabla v_t ) +
   \frac1\varepsilon ( f(u), v_t )\, dt = 0,
   \quad\forall\ v \in V^{u_0}_{[0,T]}.
\end{equation}
Here $(\cdot,\cdot)$ denotes inner product in
$L^2(\Omega)$.  By setting $v = u$ in \eqref{eq:ACweakform}, we obtain
the energy dissipation:
\begin{equation}
   \label{eq:eLaw2}
   E[u(x,T)] - E[u(x,0)] = - \int_0^T (u_t, u_t)dt \le 0.
\end{equation}
We note that taking $v\in V^{u_0}_{[0,T]}$ in \eqref{eq:ACweakform} is
equivalent to taking $v_t \in W_{[0,T]} := L^2(0,T; H_0^1(\Omega))$. We refer
the readers to \cite{feng_numerical_2003, wang_energystable_2019} for the
regularities of solutions to equation \eqref{eq:AC}. Our schemes are spectral
Galerkin methods based on weak form \eqref{eq:ACweakform}. Comparing to
existing spectral time marching method based on Petrov-Galerkin or
weighted-Galerkin approach~\cite{tang_single_2006, shen_fourierization_2007},
our approach has two obvious advantages: 1) If the continuous equation
conserves total mass then the Galerkin methods based on weak form
\eqref{eq:ACweakform} does. 2) Energy dissipation is kept like
\eqref{eq:eLaw2}. By utilizing spectral/spectral-element approximation in both
time and space, we obtain schemes with an adjustable order and good energy
stability.

The remain part of this paper is organized as follows. In Section 2, we present
the numerical schemes: one implicit scheme and one semi-implicit scheme, with
their major properties and some implementation details. In Section 3, we
present standard error estimates for the proposed schemes. Section 4 is devoted
to the superconvergence of the implicit scheme and how it can be efficiently
solved by using a Picard-like iteration similar to the semi-implicit scheme.
Section 5 presents numerical results, where energy stability and
superconvergence are numerically verified, and computational efficiency
comparison among the proposed methods and some other popular approaches are
made. In this section, we also solve the conservative Allen--Cahn equation,
which shows our approach can maintain mass conservation to machine accuracy. We
end the paper with some concluding remarks in Section 6.

\section{The energetic spectral-element time marching (ESET) schemes}

We now describe space-time Galerkin schemes based on the weak form
\eqref{eq:ACweakform}.

We first define finite dimension approximations of $V_{[s,t]}$ and
$V^{u_0}_{[s,t]}$:
\begin{equation*}
   V_{M,N}^{[s,t]} := P_{M}^0(\Omega) \otimes P_{N}\big([s,t]\big), \quad
   P_{M}^0(\Omega) := \left\{\, u(x)\in P_{M}(\Omega),\, u|_{\partial
     \Omega} =0 \,\right\},
\end{equation*}
\begin{equation*}
   V_{M,N}^{[s,t],u_0} := \left\{\, u\in V_{M,N}^{[s,t]},\ u(x,s)= u_0(x)\in
   P_{M}^0(\Omega) \,\right\},
\end{equation*}
where $P_M(\Omega)$ denotes all polynomials of degree no greater than
$M$ defined on domain $\Omega$.

\subsection{Galerkin approximation on a single element}
The Galerkin approximation for \eqref{eq:ACweakform} is defined as: Find $h\in
V^{[0,T], u_0}_{M,N}$, s.t.
\begin{equation}
   \label{eq:Galerkin}
   \int_0^T (
   h_t, v_t ) + \varepsilon( \nabla h, \nabla
   v_t ) + \frac1\varepsilon ( f(h), v_t
   )\, dt = 0, \quad\forall\, v \in V_{M, N}^{[0,T], u_0}.
\end{equation}
Here we use $h$ to denote numerical solutions.
By setting $v= h$ in \eqref{eq:Galerkin},
we immediately obtain the discrete energy law:
\begin{equation*}
   \label{eq:eLawGal}
   E[h(x,T)] - E[h(x,0)] = - \int_0^T (h_t, h_t)dt \le 0.
\end{equation*}
The Allen--Cahn equation \eqref{eq:AC} has an intrinsic time scale
$O(\varepsilon)$. For problems with long time period $T \gg \varepsilon$, a
numerical method based on \eqref{eq:Galerkin} would need to solve a nonlinear
system with a very large $N$, which is not efficient. Moreover, the solution
may not be unique. To obtain efficient methods, we consider instead a
spectral-element approach.

\subsection{Spectral-element Galerkin approximation: the implicit scheme}
We first partition the grid for $t\in[0, T]$: $0=t_0 < t_1 < \cdots < t_{K} =
T$.  For simplicity, we consider in this paper only uniform grids: $t_n = n
\tau, n=0,1,\ldots, K$, $\tau=T/K$.  Then, on each time interval
$I_n:=[t_{n-1}, t_{n}]$, we consider a sub-problem that is similar to
\eqref{eq:Galerkin}. Denote by $h^n (x, t)$ the numerical solution on interval
$I_n$, $u(x, t_n)$ the exact solution at $t_n$, $u^n(x) := h^n (x, t_n)$
the numerical solution at $t_n$, $V_{M,N}^{n,
u^{n-1}}:=V_{M,N}^{[t_{n-1}, t_n], u^{n-1}}$. Then, the weak form for the exact
solution on interval ${I}_n$ reads: Find $u(x, t) \in V_{[t_{n-1}, t_n]}^{u(x,
t_{n-1})}$, s.t.
\begin{equation} \label{eq:SEweak}
\int_{t_{n - 1}}^{t_n} ( u_t, v_t ) + \varepsilon (
\nabla u, \nabla v_t ) + \frac{1}{\varepsilon} ( f (u),
v_t )  d t = 0, \quad \forall \, v \in V_{[t_{n-1}, t_n]}^{u(x,t_{n-1})}.
\end{equation}
The corresponding Galerkin approximation is:
Find $h^n \in V_{M, N}^{n, u^{n - 1}},$ s.t.
\begin{equation}\label{eq:SEGal}
\int_{t_{n - 1}}^{t_n} ( h^n_t, v_t ) + \varepsilon (
\nabla h^n, \nabla v_t ) + \frac{1}{\varepsilon} ( f (h^n),
v_t )  d t = 0, \quad \forall\, v \in V_{M,N}^{n, u^{n-1}}.
\end{equation}

The above scheme leads to nonlinear algebraic systems to solve, since
$f(\cdot)$ is nonlinear. The advantageous aspect is that if the interval
$[t_{n-1}, t_n]$ is small enough, the system will has a unique solution.
Meanwhile, we can choose a smaller $N$ to represent the solution that satisfies
the energy dissipation property. These results are stated in the following
theorem.
\begin{theorem}\label{thm:Theorem1}
	Suppose
	\begin{equation} \label{eq:LipCond}
		\max_{u} |f'(u)| < L.
	\end{equation}
	For $\tau = |t_n - t_{n-1}| \le \varepsilon/L $, equation \eqref{eq:SEGal}
	has a unique solution $h^n$, and it satisfies the discrete energy law
	\begin{equation}
		E[h^n(x,t_n)] - E[h^n(x, t_{n-1})] = -\int_{t_{n-1}}^{t_n} (h^n_t, h^n_t) dt \le 0.
	\end{equation}
\end{theorem}
The energy stability can be obtained straightforwardly. We leave the proof of
uniqueness to Section~4.

\begin{remark}\label{rmk:1} The standard double-well potential
    \eqref{eq:double-wellF} doesn't satisfy condition \eqref{eq:LipCond} in
    Theorem \ref{thm:Theorem1}, because $|f'(u)|\rightarrow\infty$ as $|u|
    \rightarrow \infty$. But, the Allen--Cahn problem
    \eqref{eq:AC}-\eqref{eq:double-wellF} satisfies the maximum principle
    $|u(x,t)|\le 1$ if the initial value $|u(x,0)|\le 1$, so we can modify the
    double-well potential \eqref{eq:double-wellF} to have quadratic growth for
    $|u|>1$ without effecting the exact solution, such that the Lipschitz
    condition \eqref{eq:LipCond} is satisfied. Therefore, it has been a common
    practice (cf. e.g. \cite{kessler_posteriori_2004, shen_numerical_2010,
    condette_spectral_2011,wang_energystable_2019}) to consider the Allen--Cahn
    (and Cahn--Hilliard) equation with a truncated quadratic growth double-well
    potential. For example, following truncation is considered
    in~\cite{shen_efficient_2019}:
\begin{equation} \label{eq:truncatedF-M}
\tilde{F}(u) = \left\{
\begin{aligned}
&\frac{3M^2-1}{2}u^2-2M^3 u + \frac{1}{4}(3M^4+1),&\quad u> M,\\
&\frac{1}{4}(u^2-1)^2,&\quad u\in[-M,M],\\
&\frac{3M^2-1}{2}u^2+2M^3 u + \frac{1}{4}(3M^4+1),&\quad u< -M,
\end{aligned}\right.
\end{equation}
where $M \ge 1$ is a parameter. For the Cahn--Hilliard equation, one usually
choose $M>1$. For the Allen--Cahn equation, it is sufficient to choose $M=1$,
due to maximum principle. When $M=1$, \eqref{eq:truncatedF-M} reduces to
\begin{equation}\label{eq:truncatedF-M1}
\hat{F}(u) =
\begin{cases}
    (u+1)^2, & u<-1,\\
    (u^2-1)^2/4, & -1\le u\le 1,\\
    (u-1)^2, & u>1.
\end{cases}
\end{equation}
If we use $\tilde{F}(u)$ to replace $F(u)$, then \eqref{eq:LipCond} is
satisfied with $L =3M^2-1$.
\end{remark}

To solve the nonlinear system, it is usually necessary to use iterative
methods. For example, Newton's method can be applied, which has a fast
convergence rate. However, in each sub-iteration of Newton’s method, one needs
to solve large linear systems with variable coefficients, and the formation
and solution of such systems are very expensive.

We propose next a linear scheme with constant coefficients
which can be solved efficiently.

\subsection{Spectral-element Galerkin approximation: the semi-implicit linear scheme}
Now, on each sub-interval $I_n$, we consider the following semi-implicit
Galerkin approximation: Find $h^n \in V_{M, N}^{n, u^{n - 1}},$ s.t.
\begin{equation}\label{eq:SEGalLinear}
\int_{t_{n - 1}}^{t_n} ( h^n_t, v_t ) + \varepsilon (
\nabla h^n, \nabla v_t ) + \frac{S}{\varepsilon} (
    h^n, v_t ) +  \frac{1}{\varepsilon} ( \hat{f} (\hat{h}^n),
v_t )  d t = 0, \quad \forall\, v \in V_{M,N}^{n, u^{n-1}}.
\end{equation}
Here $S>0$ is a stabilization constant to ensure that $\hat{f}(u) := f(u) - S
u$ is concave within the specified region of $u$~\cite{shen_numerical_2010}.
For the Allen--Cahn equation, the maximum principle implies that $u(x, t)\in
[-1,1]$ if the initial condition $u_0(x) \in [-1, 1]$. By choosing $S=2$,
$\hat{f}$ is concave for $u\in [-1,1]$.  In \eqref{eq:SEGalLinear}, $\hat{h}^n$
is defined as the extension of $h^{n-1}$ from $I_{n-1}$ to $I_n$, for $n > 1$.
For the first step $n=1$, we can use an implicit scheme or other existing
high-order schemes to obtain $h^1$. It is easy to verify that
\eqref{eq:SEGalLinear} forms a linear positive-definite system for the unknown
degree of freedom. Note that a single stabilization term can make first order
schemes unconditionally energy stable~\cite{shen_numerical_2010}, but it
is not enough to guarantee unconditional stability for a second-order
scheme~\cite{wang_efficient_2018}. Here, we can not prove unconditionally
energy stable by adding the stabilization term, but it indeed has some
stabilization effect numerically, so we choose to keep \eqref{eq:SEGalLinear}
in the current form with $S>0$ being a tuning parameter.

\subsection{Implementation}

Now we turn to numerical implementation.
We use \eqref{eq:SEGalLinear} as an example to describe implementation
details.
We first use a time rescaling
\begin{equation*}\label{eq:tcellmap}
   t = \frac{1-\xi}{2} t_{n-1} + \frac{1+\xi}{2} t_{n},\quad \xi\in
   I:=[-1,1],
\end{equation*}
to transform \eqref{eq:SEGalLinear} on interval $I_n$ to a problem
defined on the standard reference interval $I$:
Find $h^n\in U^{h^{n}_0}_{M, N}$, such that
%\begin{subequations}
  \begin{align*}
  \label{eq:LseGalerkin}
   \int_{-1}^1 \frac2\tau  ( h^n_\xi, v_\xi ) + \varepsilon(
   \nabla h^n, \nabla v_\xi ) + \frac{S}{\varepsilon} (
   h^n, v_\xi ) + \frac1\varepsilon (
   \hat{f}(\hat{h}^n), v_\xi )\, d\xi &=0, \quad\forall\,v \in
   U_{M, N}^{h^{n}_0}, \\
   \hat{h}^n(x,\xi) = h^{n-1}(x,\xi+2),
   \qquad h^n(x,\xi=-1) = h^n_0(x) &:= h^{n-1}(x,1),
 \end{align*}
%\end{subequations}
where $U_{M,N} := V_{M,N}^{[-1,1]}$,
$U_{M, N}^{v} :=  V_{M,N}^{[-1,1], v} $.

We use following compactly-combined
Legendre bases for temporal discretization
\begin{equation}\label{eq:CompactLegendBasis}
	\phi_0(\xi)= 1, \qquad
	\phi_1(\xi)= (1+\xi)/2, \qquad
	\phi_k(\xi) = L_k(\xi) - L_{k-2}(\xi), \quad k=2,3,\ldots.
\end{equation}
We note that this kind of compactly combined Legendre bases was introduced by
Shen in~\cite{shen_efficient_1994}, where $\{\phi_k, k\ge 2\}$ are used to
treat elliptic equations with Dirichlet boundary conditions. Here we include
$\phi_0, \phi_1$ to treat initial value problem. Since
$\phi_k(-1)=\phi_k(1)=0$ for $k\ge 2$, so only $\phi_0(\xi)$ and $\phi_1(\xi)$
are used to (efficiently) interpolate the values of a solution at $\xi =-1, 1$,
which correspond to the initial value of current time interval and next time
interval. More details are given in \eqref{eq:testfunc}-\eqref{eq:LAS}.
Similar bases were used for problems with natural and other boundary conditions
\cite{wang_efficient_2018,yu_numerical_2017,yang_efficient_2018a,
wang_convergence_2018a}. We note that, the $\tau$-methods proposed by
\cite{tang_single_2002} is similar to our method applied to linear equations,
but they use different bases leading to dense matrices.

Different types of Galerkin methods can be used for spatial discretization. For
simplicity, we focus on one-dimensional and two-dimensional rectangular spatial
domains. This allow us to employ similar and highly accuracy spectral methods,
ensuring that the spatial discretization error is smaller than temporal
discretization error. We use $\psi_j(x), j=1,\ldots, M$ to denote the spatial
bases. Then the numerical solution $h^n(x, \xi)$ can be represented as:
\begin{equation}\label{eq:testfunc}
   h^{n}(x,\xi)
  =  \sum_{i=0}^{N}\sum_{j=1}^{M}\tilde{h}^{n}_{ij}\phi_i(\xi)\psi_j(x),
        \quad \xi \in [-1,1], \quad x \in \mathbb{R}^d,
\end{equation}
where $h^n(x,-1) = \sum_{j=1}^{M} \tilde{h}^n_{0j} \psi_j(x)$ is known.
Denote
\begin{align*}
    A_x &= (a^x_{ij}), &
    A_\xi&=(a^\xi_{ij}), &
    B_x &= (b^x_{ij}), &
    C_\xi &= (c^\xi_{ij}) \\
    a^x_{ij}&=(\psi'_i(x),\psi'_j(x)), &
    a^\xi_{ij}&= (\phi'_i(\xi),\phi'_j(\xi)),&
    b^x_{ij}&=(\psi_i(x),\psi_j(x)),&
    c^\xi_{ij}&= (\phi'_i(\xi),\phi_j(\xi)),
\end{align*}
with $i=1,\ldots, N, j=1,\ldots, M$. Then in each time interval $I_n$,
equation \eqref{eq:SEGalLinear} leads to a linear algebraic system
\begin{equation}\label{eq:LAS}
\frac{2}{\tau} A_\xi H^{n} B_x+\varepsilon C_\xi H^{n} A_x
+ \frac{S}{\varepsilon} C_\xi H^n B_x = \frac{1}{\varepsilon}F+ R,
\end{equation}
where $H^{n} = \bigl(\tilde{h}^{n}_{ij}\bigr)$,
$F =\bigl(\tilde{f}_{ij}\bigr)$,
$R=(r_{ij})$,
and
\begin{align*}
    \tilde{f}_{ij}
    &= -\int_{-1}^{1}( \hat{f}(\hat{h}^{n}),\psi_j(x))\phi_i'(\xi) d\xi, \\
    r_{ij}
    & = -\varepsilon  \int_{-1}^{1}
        ( \nabla h^{n}(x, -1),\nabla\psi_j(x) ) \phi_i'(\xi) d\xi
     -\frac{S}{\varepsilon} \int_{-1}^{1}
     (h^{n}(x, -1), \psi_j(x) ) \phi_i'(\xi) d\xi.
\end{align*}
Note that, $R$ is related to the initial condition $h^n(x,-1)$ on interval
$I_n$, which is known by the continuous condition $h^n(x,-1)=h^{n-1}(x, 1)$.

We have two efficient ways to solving the resulting linear system
\eqref{eq:LAS}. The first one is to use a direct sparse solver, since $A_\xi$
is diagonal, $A_x$  and $B_x$ are sparse, and $C_\xi$ is a sparse matrix which
has $2N-1$ non-zero elements in the forming form
\begin{equation*}
C_{\xi} =
\begin{pmatrix}
a_0 & a_1 &  & &   \\
-a_1 & 0 & a_2 &  &\\
 & -a_2 & 0 & \cdots &\\
 &  & \cdots & \cdots & a_{N-1}\\
 &  &   &a_{N-1}  & 0
\end{pmatrix},
\end{equation*}
where $a_i, i = 0,1,\cdots, N-1$ are non-zeros. The second way is to use the
matrix diagonalization approach~\cite{lynch_direct_1964}, which  involves
solving the generalized eigenvalue problem:
\begin{equation*}
  A_\xi E = C_\xi E\Lambda,
\end{equation*}
where $E$ is formed by generalized eigenvectors and $\Lambda$ is a diagonal
matrix composed of generalized eigenvalues. Note that the eigenvalues and
eigenvector are complex valued, since $C_\xi$ is not symmetric. It can be
proved that such a diagonalization always exist, see
\cite{kong_eigenvalue_2022} fore more details. Setting $H^{n} = EV$ in
\eqref{eq:LAS}, we get
\begin{equation*}
 \frac{2}{\tau} C_\xi E\Lambda VB_x + \varepsilon C_\xi EVA_x +\frac{S}{\varepsilon}
 C_\xi E VB_x
  = \frac{1}{\varepsilon}F + R.
\end{equation*}
Multiplying the above equation by $E^{-1}C^{-1}_\xi$, we arrive at
\begin{equation*}
 \frac{2}{\tau} \Lambda V B_x +\varepsilon V A_x + +\frac{S}{\varepsilon}
  VB_x
 = E^{-1}C^{-1}_\xi(\frac{1}{\varepsilon}F + R),
\end{equation*}
which is a series of independent systems that involves only sparse spatial
matrices, thus can be solved efficiently using fast spatial solvers. We note
that, for large $N$, there exist other type of fast solution
methods~\cite{chen_efficient_2023}.

\section{Stability and convergence analysis}

\subsection{Some preliminaries}
Denote $H_0^1(I) = \{u\in H^1(I): u(\pm 1)=0\}$, $P_N^0(I) = \{u\in P_N:u(\pm
   1)=0\}$. We first present some preliminary results.

\begin{lemma}
  \label{lem:1}Suppose that $h \in H^1 ([0, \tau])$, then
  \begin{equation}
     \int_0^{\tau} | h (t) |^2 d t \le 2 \tau | h (0) |^2 + 2 \tau^2
     \int_0^{\tau} | h_t (t) |^2 d t. \label{eq:lem1e1}
  \end{equation}
  \begin{equation}
     \int_0^{\tau} | h (t) |^2 d t \le 2 \tau | h (\tau) |^2 + 2
     \tau^2 \int_0^{\tau} | h_t (t) |^2 d t. \label{eq:lem1e2}
  \end{equation}
  If $h (0) = 0$ or $h (\tau) = 0$, the estimate can be improved as
  \begin{equation}
     \int_0^{\tau} | h (t) |^2 d t \le \tau^2 \int_0^{\tau} | h_t (t)
     |^2 d t. \label{eq:lem1e3}
  \end{equation}

  \begin{proof}
      Let $h (t) = h (0) + g (t)$, then $g (0) = 0$, and
      \begin{equation}
         \int_0^{\tau} | h (t) |^2 d t \le 2 \tau | h (0) |^2 + 2
         \int_0^{\tau} | g (t) |^2 d t \label{eq:lm1pe1} .
      \end{equation}
      For the last term in above equation, we have
      \[ \int_0^{\tau} | g (t) |^2  d t = \int_0^{\tau} \left| \int_0^t g_t
            (s) d s \right|^2 d t \le \int_0^{\tau} \left(
      \int_0^{\tau} | g_t (s) | d s \right)^2 d t \le \int_0^{\tau}
      \left( \sqrt{\tau} \| g_t \| \right)^2 d t = \tau^2 \| g_t
      \|^2, \] where H\"older inequality is used. Substituting the above
      inequality into \eqref{eq:lm1pe1}, we obtain
      \eqref{eq:lem1e1}. The inequality \eqref{eq:lem1e2} is obtained
      by symmetry. \eqref{eq:lem1e3} is obtained by
      noticing $h (t) = g (t)$ with $g (0) = 0$ or $g (\tau) = 0$.
  \end{proof}
\end{lemma}

\begin{lemma}
   \label{lem:Lextra}{} Denote by $L_k (x)$ the Legendre polynomial of degree
    $k$. Then there exist constants $1 \le c_k < 2^{2k-1}3^{2k+1}$, such
    that
   \begin{equation}
     \int_1^3 | L_k (x) |^2  dx\
        = \ c_k \int_{- 1}^1 | L_k (x) |^2 dx,
        \quad k=0, 1, \ldots.
   \end{equation}
    In particular
   \begin{equation*}
        c_0 = 1, \quad c_1 = 13, \quad c_2 = 241, \quad c_3 = 5629.
   \end{equation*}
   \begin{proof}
      Since $L_k (x), k = 0, 1, \ldots$ are orthogonal polynomials defined on
      the interval $[- 1, 1]$, so their roots all lie in $[- 1,
        1]$. Furthermore, by $L_k (1) = 1$, we know that $L_k (x) > 0$
      for $x > 1$. By $L_0=1, L_1(x)=x$, and the three term recurrence
     \[ L_{k + 1} (x) = \frac{2 k + 1}{k + 1} xL_k (x) - \frac{k}{k + 1} L_{k
                - 1} (x), \quad k \ge 1,\]
    we have $L_{k + 1} (x) < 2 xL_k (x)$, for
     $x \in [1, 3]$. By $L_0 (x) \equiv 1$ and mathematical induction,
     we get $L_k (x) \le (2 x)^k$, for
     $x \in [1, 3]$.  Therefore,
     \[ \int_1^3 L_{k}^2 (x) dx \le \int_1^3 2^{2 k} x^{2 k}
     dx = \frac{2^{2 k}}{2 k + 1} (3^{2 k + 1} - 1) < \frac{2^{2
         k}}{2 k + 1} 3^{2 k + 1},\quad k>0.\] Since $\int_{- 1}^1 L_k^2 (x)
     dx = \frac{2}{2 k + 1}$, so $c_k < 2^{2 k - 1} 3^{2 k +
       1}$, for $k>0$. The exact values for small $k$'s can be obtained by direct
     calculations. Obvious $c_0 = 1$. For $k = 1$, $L_1 (x) = x$,
     \[ \int_1^3 x^2  dx = \frac{1}{3} x^3 \mid_1^3 = \frac{26}{3}, \qquad
       \text{which leads to} \quad c_1 = \frac{26}{3} / \frac{2}{3} = 13. \] For
       $k = 2$, $L_2 (x) = \frac{1}{2} (3 x^2 - 1)$
     \[ \int_1^3 \frac{1}{4} (9 x^4 - 6 x^2 + 1)  dx = (\frac{9}{20} x^5 -
       \frac{1}{2} x^3 + \frac{x}{4} ){\mid}_1^3 = \frac{482}{5}, \quad
       \text{which leads to}  \quad c_2 = \frac{482}{5} / \frac{2}{5} = 241. \]
       For $k = 3$, $L_3 (x) = \frac{1}{2} (5 x^3 - 3 x)$
       \[ \int_1^3 \frac{1}{4} (5 x^3 - 3 x)^2  dx = \frac{11258}{7}, \quad
            \text{which leads to}  \quad c_3 = \frac{11258}{7} / \frac{2}{7} =
          5629. \]
   \end{proof}
\end{lemma}

   \begin{lemma}
    \label{lem:3}Suppose that $h \in P_N ([0, 2\tau])$, then
        \begin{equation}
            \label{eq:lem3}
            \int_{\tau}^{2 \tau} | h (t) | dt
            \le C_N \sqrt{\tau}
                \left( \int_0^{\tau} | h (t) |^2 dt \right)^{1/2},
        \end{equation}
        where $C_N = \sqrt{\sum_{k = 0}^N c_k}$. In particular $C_0 = 1$, $C_1 =
        \sqrt{14} \cong 3.7, C_2 = \sqrt{255} \sim 16.0, C_3 = \sqrt{5884} \sim
        76.7$.

        \begin{proof}
            First define
            \[ \varphi_k (t) = \sqrt{\frac{2 k + 1}{\tau}} L_k \left( \frac{2 t}{\tau}
            - 1 \right), \]
            where \ $L_k (x)$ is Legendre polynomial of degree $k$. It is easy to
            verify that $\{ \varphi_k (t) \}$ compose of a set of orthonormal bases in
            $L^2 ([0, \tau])$:
            \[ \int_0^{\tau} \varphi_k (t) \varphi_j (t)  dt = \int_{- 1}^1
            \frac{\sqrt{(2k+1)(2j+1)}}{\tau} L_k (\xi) L_j (\xi)
            \frac{\tau}{2}  d\xi
            %= \frac{2 k + 1}{2} \int_{- 1}^1 L_k (\xi) L_j (\xi)  d\xi
            = \delta_{j k} . \]
            Then, by using Lemma \ref{lem:Lextra}, we have
            \begin{equation}
            \int_{\tau}^{2 \tau} | \varphi_k (t) |^2  dt = \frac{2 k + 1}{2}
            \int_1^3 | L_k (\xi) |^2  d\xi = c_k \frac{2 k + 1}{2}
            \int_{- 1}^1 | L_k (\xi) |^2  d\xi = c_k .
            \end{equation}
            By expanding $h (t)$ as
            $ h (t) = \sum_{k = 0}^N b_k \varphi_k (t)$,
            we get
            \[ \int_0^{\tau} | h (t) |^2  dt = \sum_{k = 0}^N b_k^2, \]
            and
            \begin{align*}
                \int_{\tau}^{2 \tau} | h (t) |   dt
                & \le  \sum_{k = 0}^N
                |b_k| \int_{\tau}^{2 \tau} | \varphi_k |  dt \le  \sum_{k = 0}^N
                |b_k| \sqrt{\tau}  \left( \int_{\tau}^{2 \tau} | \varphi_k |^{^2}
                dt \right)^{1 / 2}\\
                & = \sum_{k = 0}^N |b_k|  \sqrt{\tau} \sqrt{c_k} \le
                \sqrt{\tau \sum_{k=0}^N c_k} \left( \int_0^{\tau} | h (t) |^2  dt
                \right)^{1 / 2}\\
                & =  C_N \sqrt{\tau} \left( \int_0^{\tau} | h (t) |^2  dt
                \right)^{1 / 2} .
            \end{align*}
            The proof is complete.
        \end{proof}
    \end{lemma}

    Next we present several lemmas regarding projection error.
    Let $w^{\alpha,\beta}(t) := (1-t)^{\alpha} (1+t)^{\beta} $ with
    $\alpha,\beta > -1$. Introduce the non-uniformly Jacobi-weighted Sobolev
    space:
    \begin{equation}
    	B^m_{\alpha,\beta}(I)
        = \big\{\, u: \partial^k_t u \in L^2_{w^{\alpha+k,\beta+k}}(I),\ 0\le k\le m \,\big\},
            \ m\in \mathbb{N},
    \end{equation}
    equipped with the inner product, norm and semi-norm
    \begin{equation*}
    	(u,v)_{B^m_{\alpha,\beta}}
        = \sum_{k=0}^{m}(\partial^k_t u,\partial^k_t v)_{w^{\alpha+k,\beta+k}},
    \end{equation*}
    \begin{equation*}
    	\lVert u \rVert_{B^m_{\alpha,\beta}}
        = (u,u)^{1/2}_{B^m_{\alpha,\beta}},\ |u|_{B^m_{\alpha,\beta}}
        = \lVert \partial^m_t u\rVert_{w^{\alpha+m,\beta+m}}.
    \end{equation*}

	For given temporal or spatial domain $\Omega$, we define
	\begin{equation*}
	    a(u,v) = \int_\Omega \nabla u(\xi) \cdot \nabla v(\xi) d\xi.
	\end{equation*}
	Consider the  orthogonal projection $\Pi^{0}_N: H_0^1(\Omega)\mapsto
	P_N^0(\Omega)$, defined by
	\begin{equation}
	a(\Pi^{0}_N u-u, v) = 0,\quad \forall v\in P_N^0(\Omega).
	\end{equation}
	For the one-dimensional case, we have the following lemma.  % [Theorem 3.38]{shen_spectral_2011}
	\begin{lemma}{\cite{shen_spectral_2011, canuto_spectral_2007}}\label{Lemma:pi0}
		If $u\in H_0^1(I)$ and $\partial_\xi u\in B^{m-1}_{0,0}(I)$ ,
		then for $1\le m\le N+1$ and $ \mu = 0,1$,
		\begin{equation}
		\lVert\Pi^{0}_N u-u\rVert_\mu\le
		c\sqrt{\frac{(N-m+1)!}{N!}} (N+m)^{\mu-(m+1)/2}
		\lVert\partial_\xi^m u\rVert_{w^{m-1,m-1}},
		\end{equation}
		where  $c$ is a positive constant independent of $m,N$ and $u$.
	\end{lemma}

	Define linear interpolation operator $\mathcal{I}^n_1: C^0(I_n) \mapsto
	P_1(I_n)$ as
	\begin{equation}
	    \mathcal{I}_1^n u(t)
        = \frac{t_n-t}{\tau}u(t_{n-1}) + \frac{t - t_{n-1} }{\tau}u(t_n).
	\end{equation}
	Define projection $\Pi_N^{0,n} :H^1_0(I_n) \mapsto P_N^0(I_n)$ as
	\begin{equation}
	    \int_{I_n} (u- \Pi_N^{0,n} u)'(t)\, v'(t)\, dt
        =0, \quad  \forall\, v\in P_N^0(I_n).
	\end{equation}
	Then, define operator $\Pi_N^{n} :H^1(I_n) \mapsto P_N(I_n)$ as
	\begin{equation}
	 \Pi_N^n u(t) = \mathcal{I}^n_1 u(t) + \Pi_N^{0,n} (u - \mathcal{I}^n_1 u).
	\end{equation}

	By using Lemma \ref{Lemma:pi0}, the following result can be obtained.
	\begin{lemma}\label{Lemma:PiNn}
		If $u\in H^1(I_n)$ and $\partial_t u\in B^{m-1}_{0,0}(I_n)$,
		then for $\mu = 0,1$, we have
		\begin{equation}
		\begin{aligned}
		\lVert \partial_t^\mu (\Pi^{n}_N u-u)\rVert_{L^2(I_n)}
		\le 	c\left(\frac{\tau}{2}\right)^{m-\mu}\sqrt{\frac{(N-m+1)!}{N!}}
                (N+m)^{\mu-(m+1)/2}
		        \lVert\partial_t^m u\rVert_{L^2_{w^{m-1,m-1}}(I_n)},
		\end{aligned}
		\end{equation}
		for  $1\le m\le N+1$,
		where  $c$ is a positive constant independent of $m,N$ and $u$.
		In particular, when $m=N+1$, we have
		\begin{equation}\label{ineq:pi}
		\lVert \partial_t^\mu (\Pi^{n}_N u-u)\rVert_{L^2(I_n)}
		\le  c\left(\frac{\tau}{2}\right)^{N+1-\mu}
		    (2\pi N)^{-1/4}\biggl(\frac{\sqrt{e/2}}{N}\biggr)^N
		    \left(\frac{1}{2N + 1}\right)^{1-\mu}
		    \lVert\partial_t^{N+1} u\rVert_{L^2_{w^{N,N}}(I_n)}.
		\end{equation}
	\end{lemma}

  	In space, we consider similar projection $\Pi_M^0: H_0^1(\Omega)\rightarrow
  	P_M^0(\Omega)$. We assume the solution is smooth and a large enough $M$ is
  	used, such that the spatial projection error is small comparing to temporal
  	discretization error:
  	\begin{equation}\label{eq:sperr-assump}
  		\lVert\nabla^\mu (\Pi^{0}_M u - u)\rVert \le c M^{-r}
        \le c \tau^{2N}, \quad \mu=0, 1.
  	\end{equation}

\subsection{Energy stability of the semi-implicit linear ESET method}
We first introduce following shorthand notations:
\begin{equation}
(u, v)_n := \int_{t_{n-1}}^{t_{n}} \int_{\Omega}uv\ dxdt, \quad
\ \lVert u\rVert_n^2:=\int_{t_{n-1}}^{t_{n}}\int_{\Omega}u^2\ dxdt.
\end{equation}
We use $u(t)$ in places of $u(x,t)$ for notation simplicity.
%In the analysis, we use assume that
%\begin{equation} \label{eq:fLip}
%\max_u |f'(u)| \le L,
%\end{equation}
%which is reasonable, since the solution of Allen--Cahn equation satisfies
%maximum principle, we can modify $f(u)$ for $x\notin [-1,1]$ without affect
%the solution to make \eqref{eq:fLip} true.
Our first result regarding the energy stability of scheme
\eqref{eq:SEGalLinear} is given in following theorem.
\begin{theorem}
	Under the condition
	\begin{equation} \label{eq:stepcond}
	\tau \frac{L}{\varepsilon} \le \frac{\sqrt{14}}{4\sqrt{2+C_{N-1}^2}},
	\end{equation}
    the following energy dissipation law
	\begin{equation} \label{eq:disEnergyLaw}
	\hat{E}^{n+1}\le\hat{E}^n, \ \forall\, n\ge 1,
	\end{equation}
	holds for scheme \eqref{eq:SEGalLinear}, where
	\begin{equation}
	    \hat{E}^{n}
        = E(h^n(t_n))
            + \frac{\tau^2L^2C_{N-1}^2}{2\varepsilon^2}\lVert h_t^n \rVert_{n}^2.
	\end{equation}
\end{theorem}
\begin{proof}
	We only consider the case $S=0$. The analysis to case $S\neq 0$ is similar.
    By taking % difference of equation \eqref{eq:SEweak} and
    $v_t=h^n_t$ in \eqref{eq:SEGalLinear}, we obtain
	\begin{equation}\label{eq:esetLinErr}
	\lVert  h_t^n \rVert_{{n}}^2 + (\frac{d}{dt}E(h^n),1)_{{n}}
	= \frac{1}{\varepsilon}(f(h^n)-f(\hat{h}^n),  h^n_t)_{n}
	\le \frac{L}{\varepsilon}(| h^n-\hat{h}^n|, |h_t^n|)_{{n}}.
	\end{equation}
	Here the that $f(u)=\hat{f}(u)$ when $S=0$ is used. Since
	$h^n(t_{n-1}) = \hat{h}^n(t_{n-1})$, we have
	\begin{equation*}
	| (h^n-\hat{h}^n)(t)| \le \int_{t_{n-1}}^{t}| h^n_t(s)-\hat{h}^n_t(s) | ds
	\le \int_{t_{n-1}}^{t}| h^n_t(s)| + |\hat{h}^n_t(s) | ds, \quad t\in [t_{n-1},t_{n}].
	\end{equation*}
	Plugging the above result into \eqref{eq:esetLinErr}, then for the first
	term on the right hand side, we have
	\begin{equation*}
	(\int_{t_{n-1}}^{t} |h^n_t| ds,| h^n_t|)_n = \frac12 (\int_{t_{n-1}}^{t_{n}}| h^n_t|\ ds,| h^n_t|)_n
	\le \frac12(\sqrt{\tau}\lVert h^n_t\rVert_n,| h^n_t|)_n
	\le \tau^2\lVert h^n_t\rVert_n^2 + \frac{1}{16} \lVert h^n_t\rVert_n^2,
	\end{equation*}
	where H\"older and Cauchy inequalities are used. For the second term, we do
	similar estimate and using Lemma \ref{lem:3} to get
	\begin{equation}
	\begin{aligned}
	    \bigl(\int_{t_{n-1}}^{t} |\hat{h}^n_t| ds,|h^n_t|\bigr)_n
        \le (C_{N-1}\sqrt{\tau}\lVert h^{n-1}_t\rVert_{n-1},| h^n_t|)_n
	    \le \frac{\tau^2 C_{N-1}^2}{2} \lVert h^{n-1}_t\rVert_{n-1}^2
            + \frac{1}{2}\lVert h^n_t\rVert_n^2,
	\end{aligned}
	\end{equation}
	Putting these results into \eqref{eq:esetLinErr}, we obtain the discrete
	energy dissipation relation
	\begin{equation*}
		E(h^n(t_n)) + \frac{\tau^2 C_{N-1}^2 L^2}{2\varepsilon^2} \| h^{n}_t \|_{n}^2
        \le  E(h^n(t_{n-1})) + \frac{\tau^2 C_{N-1}^2 L^2}{2\varepsilon^2} \| h^{n-1}_t \|_{n-1}^2
		 - \left( \frac{7}{16} - \frac{\tau^2 L^2}{\varepsilon^2} \frac{C_{N-1}^2+2}{2} \right)  \| h^n_t \|_n^2.
	\end{equation*}
    We obtain desired energy dissipation under condition \eqref{eq:stepcond}.
\end{proof}

\subsection{Convergence analysis}

\subsubsection{A standard error estimate for the implicit ESET scheme}

Let $h^n$, $u(t)$ be the solution of equation \eqref{eq:SEGal} and
\eqref{eq:SEweak}, correspondingly. Denote by
\[ e ^n (t) = h^n (t) - u (t) = \rho^n + \pi^n, \quad \rho^n (t) = h^n (t)
-\Pi^n u (t), \quad \pi^n  (t) =\Pi^n u (t) - u (t), \] where $\Pi^n u :=
\Pi_N^n \Pi_M^0 u$  is a spatial-temporal projection operator. By Lemma
\ref{Lemma:PiNn} and assumption \eqref{eq:sperr-assump}, we have
\begin{equation}
\| \nabla^\mu \pi^n \|
\le \| \nabla^\mu (u - \Pi_M^0 u) \| + \| \nabla^\mu (\Pi_M^0 u - \Pi_N^n \Pi_M^0 u)\|
\le  c M^{-r} + c \tau^{N+1},\quad  \mu=0,1,
\end{equation}
where $c$ is a general constant. Since $\|\nabla^\mu  e^n\| = \|\nabla^\mu
\rho^n\| + \|\nabla^\mu  \pi^n\|$, to get a $H_0^1(\Omega)$ upper bound of the
numerical error, we only need to estimate $\|\nabla^\mu  \rho^n\|$, for which
we have the following result.
\begin{theorem}\label{th:inpeset}
	Following error estimates hold for the implicit ESET scheme \eqref{eq:SEGal}:
\begin{equation}\label{eq:ErrImESST}
	\| \rho^n (t_{n}) \|^2 + \frac{\varepsilon^2}{2 L} \| \nabla \rho^n (t_{n}) \|^2
    + \sum_{k = 1}^n D_{\tau} \| \rho_t^k \|_k^2 \le \sum_{k = 1}^n e^{D_2 k
		\tau} E_p^k (\tau), \quad \forall\, \tau \le \frac{\varepsilon}{4
		\sqrt{2} L} .
\end{equation}
where $D_{\tau} = \frac{\varepsilon}{4L} - \frac{8 \tau^2 L}{\varepsilon}$,
$D_2 = \frac{8 L }{\varepsilon} $, and
\begin{equation}\label{eq:projerrineq}
    E^n_p(\tau)
    =  6\| \partial_t (\Pi_M^0-I_{id}) \Pi_N^n u \|_n^2
     + 6 \varepsilon^2 \| \nabla^2 (\Pi_N^n-I_{id}) \Pi_M^0 u \|_n^2
     + \frac{6 L^2}{\varepsilon^2} \| \pi^n \|_n^2,
\end{equation}
where $I_{id}$ stands for identity operator.
By Lemma \ref{Lemma:PiNn} and assumption on spatial projection error
\eqref{eq:sperr-assump}, we have $E_p^n(\tau) < \mathcal{O}(\tau M^{-2r} +
\tau^{2N+3})$. So the $H_0^1(\Omega)$ numerical error at a time $t$ is of order
$\mathcal{O}(M^{-r}+\tau^{N+1})$.
\end{theorem}
\begin{proof}
By taking the difference of \eqref{eq:SEweak} and \eqref{eq:SEGal}, we get
\[ ( h^n_t - u_t, v_t )_n + \varepsilon
( \nabla (h^n - u), \nabla v_t )_n + \frac{1}{\varepsilon} (f
(h^n) - f (u), v_t)_n = 0, \]
which leads to
\[  ( \rho^n_t + \pi^n_t, v_t )_n + \varepsilon ( \nabla (\rho^n  + \pi^n),
\nabla v_t )_n + \frac{1}{\varepsilon} (f' (\zeta^n) (\rho^n  + \pi^n ), v_t)_n
= 0, \] where $\zeta^n = \zeta^n (x, t)$ is some value between $u (x, t)$ and
$h^n (x, t)$. By taking $v = \rho^n$, we get
\begin{equation}
( \rho^n_t + \pi^n_t, \rho^n_t )_n +
\varepsilon ( \nabla (\rho^n  + \pi^n ), \nabla \rho^n_t )_n +
\frac{1}{\varepsilon} (f' (\zeta^n) (\rho^n  + \pi^n ), \rho^n_t)_n = 0.
\end{equation}
Rearranging the above equation yields
\begin{equation}\label{eq:errR1R2R3}
    (\rho_t^n, \rho_t^n)_n
        + \frac{\varepsilon}{2} \| \nabla \rho^n (t_n) \|^2
        - \frac{\varepsilon}{2} \| \nabla \rho^n (t_{n-1}) \|^2
    = R_1 + R_2 + R_3,
\end{equation}
where
Using the definition of projection and Cauchy-Schwartz inequality, we obtain
\begin{equation}
| R_1 | = |(\partial_t(\Pi_M^0-I_{id}) \Pi_N^n u, \rho_t^n)_n|
\le \frac{1}{2 \eta_1} \| \partial_t \pi^n_M \|_n^2 + \frac{\eta_1}{2} \|
\rho_t^n \|_n^2, \quad 0 < \eta_1, \label{eq:R1}
\end{equation}
\begin{equation}
    \label{eq:R2}
    | R_2 |
    = |- \varepsilon (\nabla ( (\Pi_N^n-I_{id}) \Pi_M^0 u), \nabla \rho_t^n)_n |
    \le \frac{\varepsilon^2}{2 \eta_2} \| \nabla^2\pi_N^n \|_n^2
        + \frac{\eta_2}{2} \| \rho_t^n \|_n^2, \quad 0 < \eta_2,
\end{equation}
\begin{equation}
| R_3 | \le \frac{L^2}{\varepsilon^2} \frac{1}{2 \eta_3} \| \rho^n
\|_n^2 + \frac{\eta_3}{2} \| \rho_t^n \|_n^2 + \frac{L^2}{\varepsilon^2}
\frac{1}{2 \eta_4} \| \pi^n \|_n^2 + \frac{\eta_4}{2} \| \rho_t^n \|_n^2,
\quad 0 < \eta_3, \eta_4, \label{eq:rnorig}
\end{equation}
where $\pi_M^n = (\Pi_M^0-I_{id}) \Pi_N^n u$, $\pi_N^n = (\Pi_N^n-I_{id})
\Pi_M^0 u$. To control the $\| \rho^n \|_n^2$ term, we use Lemma \ref{lem:1} to
get
\begin{equation}\label{eq:rhon}
\| \rho^n \|_n^2 \le 2 \tau  \| \rho^n  (t_{n - 1}) \|_n^2 + 2 \tau^2
\| \rho^n_t \|_n^2 .
\end{equation}
Combining the estimates \eqref{eq:errR1R2R3}-\eqref{eq:rhon}, we have
\begin{align*}
	& (\rho_t^n, \rho_t^n)_n + \frac{\varepsilon}{2} \| \nabla \rho^n (t_n)
	\|^2 - \frac{\varepsilon}{2} \| \nabla \rho^n (t_{n - 1}) \|^2\\
	& \le \frac{\eta_1 + \eta_2 + \eta_3 + \eta_4}{2} \| \rho_t^n
	\|_n^2 + \frac{1}{2 \eta_1} \| \pi^n_t \|_n^2 + \frac{\varepsilon^2}{2
		\eta_2} \| \nabla^2 \pi^n_N \|_n^2 + \frac{L^2}{\varepsilon^2} \frac{1}{2
		\eta_4} \| \pi^n \|_n^2 + {\frac{L^2}{\varepsilon^2}
		\frac{1}{2 \eta_3} \| \rho^n \|_n^2}\\
	& \le  \left( \frac{\eta_1 + \eta_2 + \eta_3 + \eta_4}{2} +
	\frac{L^2 \tau^2}{\eta_3 \varepsilon^2} \right) \| \rho_t^n \|_n^2 +
	\frac{L^2 \tau}{\eta_3 \varepsilon^2} \| \rho^n (t_{n - 1}) \| ^2 \\
    &\qquad\hspace{5cm} +
	\frac{1}{2 \eta_1} \| \pi^n_t \|_n^2 + \frac{\varepsilon^2}{2 \eta_2} \|
	\nabla^2 \pi^n_N \|_n^2 + \frac{L^2}{2 \eta_4 \varepsilon^2} \| \pi^n \|_n^2 .
\end{align*}
Setting $\eta_3 = 1 / 4$, $\eta_1 = \eta_2 =  \eta_4 = 1 / 12$, we obtain
\begin{equation}
	\begin{split}
		\left( \frac{3}{4} - \frac{4 \tau^2 L^2}{\varepsilon^2} \right) \| \rho_t^n
		\|_n^2 + \frac{\varepsilon}{2} \| \nabla \rho^n (t_n) \|^2
		& - \frac{\varepsilon}{2} \| \nabla \rho^n (t_{n - 1}) \|^2 \le \frac{4 \tau
			L^2}{\varepsilon^2} \| \rho^n (t_{n - 1}) \| ^2 \\
		&\quad + 6 \| \pi^n_t \|_n^2 + 6
		\varepsilon^2 \| \nabla^2 \pi^n_N \|_n^2 + \frac{6 L^2}{\varepsilon^2} \|
		\pi^n \|_n^2 . \label{eq:e14}
	\end{split}
\end{equation}
Notice that
\begin{equation}
\begin{aligned}
\| \rho^n (t_n) \|^2 - \| \rho^n (t_{n - 1}) \|^2 &  = \int_{t_{n -
        1}}^{t_{n}} 2 (\rho^n, \rho^n_t) dt
 \le \frac{\varepsilon}{2 L} \| \rho^n_t \|_n^2 + \frac{2
    L}{\varepsilon} \| \rho^n \|_n^2\\
& \le \frac{4 L \tau }{\varepsilon} \| \rho^n (t_{n - 1}) \| ^2 +
\frac{\varepsilon}{L} \left( \frac{1}{2} + \frac{4 \tau^2
    L^2}{\varepsilon^2} \right) \| \rho^n_t \|_n^2.
\end{aligned} \label{eq:15}
\end{equation}
Summing up estimate \eqref{eq:15} with \eqref{eq:e14} multiplied by
$\frac{\varepsilon}{L}$, we get
\[ \| \rho^n (t_n) \|^2 - \| \rho^n (t_{n - 1}) \|^2 + D_{\tau} \| \rho_t^n
\|_n^2 + \frac{\varepsilon^2}{2 L} \| \nabla \rho^n (t_n) \|^2 -
\frac{\varepsilon^2}{2 L} \| \nabla \rho^n (t_{n - 1}) \|^2 \le D_2 \tau \|
\rho^n (t_{n - 1}) \| ^2 + E_p^n (\tau). \] By a discrete Gr\"onwall's
inequality and the fact $\rho^1 (t_0) = 0$, we obtain the desired estimate
\eqref{eq:ErrImESST}.
\end{proof}

\subsubsection{Error estimate for the semi-implicit linear ESET scheme}

Let $h^n$, $u(t)$ be the solution of equation \eqref{eq:SEGalLinear} and
\eqref{eq:SEweak}, correspondingly. We again do following error splitting
\[ e ^n (t) = h^n (t) - u (t) = \rho^n + \pi^n, \quad \rho^n (t) = h^n (t)
-\Pi^n u (t), \quad \pi^n  (t) =\Pi^n u (t) - u (t). \] The a priori projection
error $\pi^n$ is the same as in the implicit scheme, so we only need to
estimate the $\rho^n$ term.

The challenging part for the semi-implicit ESET scheme is the extrapolation
term $\hat{h}^n$. To make it easier, we introduce a mixed-type projection: we
still use the spatial project $\Pi_M^0$ as in the fully implicit ESET scheme,
but for temporal ``projection'', we use finite-term Taylor series. More
precisely, we define $\hat\Pi^n = \mathcal{T}^n_N \Pi_M^0$, where
$\mathcal{T}_N^n : C^{N+1} \mapsto P_N$ is defined as
\begin{equation}
    (\mathcal{T}_N^n u)(t)
    := \sum_{k=0}^N \frac{1}{k!} u^{(k)}(t_{n-1}) (t-t_{n-1})^k,
        \quad \forall\, t\in I_n\cup I_{n-1}.
\end{equation}
It is easy to show that
\begin{equation} \label{eq:TaylorErr}
    \|u - \mathcal{T}_N^n u \|_\alpha
    \le \frac{\tau^{N+1}}{(N+1)!} \| u^{(N+1)} \|_\alpha,
        \quad  \alpha = n, n-1.
\end{equation}
By using $\hat{\Pi}^n$ and a similar procedure for implicit scheme, we can get
the following result.

\begin{theorem}
    The following error estimates hold for the semi-implicit linear ESET scheme
    \eqref{eq:SEGalLinear}:
    \begin{equation}\label{eq:ImexESSTerr}
    \begin{split}
    \lVert \rho^n( t_{n})\rVert^2 + & \frac{\varepsilon^2}{2L}\lVert\nabla \rho^{n}( t_{n}) \rVert^2
    + \sum_{k=1}^{n}D_\tau \lVert\rho_t^k\rVert_k^2
    + D_3\tau \lVert\rho_t^{n}\rVert_{n}^2
    \le \sum_{k=1}^{n}e^{D_2 k\tau} \hat{E}^k_p(\tau) \\
    &+ \lVert \rho^1(t_{1})\rVert^2
    + \frac{\varepsilon^2}{2L}\lVert\nabla \rho^{1}(t_1) \rVert^2
    +  D_3\tau \lVert\rho_t^{1}\rVert_{1}^2,
    \quad
    \forall\,\tau\le \frac{\varepsilon}{4L\sqrt{6C_{N-1}^2+1}},\, n\ge 2,
    \end{split}
    \end{equation}
        where $D_\tau = \frac{\varepsilon}{L}(\frac{1}{4} - \frac{4\tau^2L^2(1+6C_{N-1}^2)}{\varepsilon^2} )$,
    $D_2 = \frac{8L}{\varepsilon}$,
    $D_3 = \frac{24 L\tau C_{N-1}^2}{\varepsilon}$,  and
    \begin{align}
    \hat{E}^k_p =  6 \| \partial_t \pi^n_M \|_n^2
    +  6{\varepsilon^2} \| \nabla^2\pi_N^n \|_n^2
    &+ \frac{4L^2\tau^2}{\veps^2} \| {\pi}^n(t_{n-1}) \|^2 \nn\\
    &+ \frac{12 L^2\tau^2 C_{N-1}^2}{\veps^2} (\| \hat{\pi}^n_t \|_{n-1}^2
    +  2\| \partial_t (\hat{\Pi}^n u - {\Pi}^{n-1} u) \|^2_{n-1}).
    \label{eq:imexEpk}
    \end{align}
\end{theorem}
\begin{proof}
    By taking the difference of (\ref{eq:SEweak}) and \eqref{eq:SEGalLinear},
    and taking $v = \rho^n_t$, we obtain
    \begin{equation*}
    (\rho_t^n+\pi_t^n,\rho^n_t)_n + \varepsilon(\nabla\rho^n,\nabla \rho^n_t)_n
    + \varepsilon(\nabla\pi^n,\nabla\rho^n_t)_n=
    -\frac{1}{\varepsilon}(f(u)-f(\hat{h}^n),\rho^n_t)_n.
    \end{equation*}
    Rearranging the above equation to obtain
    \begin{equation}\label{eq:ImexErrEq}
    \lVert\rho^n_t\rVert_n^2 + \frac{\varepsilon}{2}\lVert\nabla \rho^n( t_{n}) \rVert^2
    - \frac{\varepsilon}{2}\lVert\nabla \rho^{n}(t_{n-1}) \rVert^2
    = R_1 + R_2 +R_3,
    \end{equation}
    where $R_1 = - (\pi^n_t,\rho^n_t)_n$, $R_2 =
    -\varepsilon(\nabla\pi^n,\nabla\rho^n_t)_n$. There are same as in the
    implicit ESET scheme, can be handled similarly. The $R_3$ term is
    $$R_3 = -\frac{1}{\varepsilon}(f'(\zeta^n)(\hat{h}^n - u),\rho^n_t)_n.$$
    For $\hat{ h}^n-u$  we use following splitting
    \begin{equation*}
        \hat{ h}^n-u = \hat{\rho}^n + \hat{\pi}^n,
         \quad \hat{\rho}^n (t) = \hat{h}^n (t) -\hat{\Pi}^n u (t),
         \quad \hat{\pi}^n (t) =\hat{\Pi}^n u (t) - u (t).
    \end{equation*}
    Then we have
    \begin{align}
        |R_3| & \le \frac{L}{\varepsilon} (|\hat{\rho}^n + \hat{\pi}^n|, |\rho^n_t|)_n \nn \\
        & \le \frac{L}{\varepsilon} \bigl( |\hat{\rho}^n(t_{n-1}) + \hat{\pi}^n(t_{n-1})|
         + \int_{t_{n-1}}^t |\hat{\rho}_t^n(s) + \hat{\pi}_t^n(s)| ds  , |\rho^n_t|\bigr)_n \nn \\
        & \le \frac{L^2\tau}{\veps^2 \eta_5} (\| {\rho}^n(t_{n-1})\|^2 + \| {\pi}^n(t_{n-1}) \|^2 )
        + \frac{L^2\tau^2 C_{N-1}^2}{\veps^2 \eta_6} (\| \hat{\rho}^n_t \|^2_{n-1}
        + \| \hat{\pi}^n_t \|_{n-1}^2 ) + \frac{\eta_5 + \eta_6}{2} \| \rho^n_t\|_n^2, \label{eq:ImexR3}
    \end{align}
    where the fact $\rho^n(t_{n-1})=\hat{\rho}^n(t_{n-1})$,
    $\pi^n(t_{n-1})=\hat{\pi}^n(t_{n-1})$, and Lemma \ref{lem:3} are used. The
    term $\| \hat{\rho}^n_t \|^2_{n-1}$ can be further estimated by
    \begin{align} \label{eq:Imexhatr}
        \| \hat{\rho}^n_t \|^2_{n-1} = \| \partial_t ({h}^{n-1} - \hat{\Pi}^n u) \|^2_{n-1}
        & \le 2\| \partial_t ({h}^{n-1} - {\Pi}^{n-1} u) \|^2_{n-1}
        + 2\| \partial_t (\hat{\Pi}^n u - {\Pi}^{n-1} u) \|^2_{n-1} \nn \\
        &= 2\| \rho^{n-1}_t \|^2_{n-1}
        + 2\| \partial_t (\hat{\Pi}^n u - {\Pi}^{n-1} u) \|^2_{n-1}
    \end{align}
    Combining the estimates \eqref{eq:R1}, \eqref{eq:R2}, \eqref{eq:ImexErrEq},
    \eqref{eq:ImexR3} and \eqref{eq:Imexhatr}, then taking $\eta_5 = 1 / 4$,
    $\eta_1 = \eta_2 =  \eta_6 = 1 / 12$, we obtain
    \begin{align}
    \frac34 \| \rho_t^n \|_n^2 &+ \frac{\varepsilon}{2} \| \nabla \rho^n (t_n) \|^2 -
    \frac{\varepsilon}{2} \| \nabla \rho^n (t_{n-1}) \|^2 \nn \\
    & \le \frac{4L^2\tau}{\veps^2} \| {\rho}^n(t_{n-1})\|^2
    + \frac{24 L^2\tau^2 C_{N-1}^2}{\veps^2} \| \rho^{n-1}_t \|^2_{n-1}
    + 6 \| \partial_t \pi^n_M \|_n^2
    +  6{\varepsilon^2} \| \nabla^2\pi_N^n \|_n^2  \nn \\
    &\quad + \frac{4L^2\tau}{\veps^2} \| {\pi}^n(t_{n-1}) \|^2
    + \frac{12 L^2\tau^2 C_{N-1}^2}{\veps^2} (\| \hat{\pi}^n_t \|_{n-1}^2
    +  2\| \partial_t (\hat{\Pi}^n u - {\Pi}^{n-1} u) \|^2_{n-1}).
    \end{align}
    Multiplying the above equation by $\frac{\veps}{L}$, then summing the
    results with \eqref{eq:15}, we reach to
    \begin{align}
     \| \rho^n (t_n) \|^2 -& \| \rho^n (t_{n - 1}) \|^2
     + D_{\tau} \| \rho_t^n \|_n^2
     + D_3 \tau \| \rho_t^{n} \|_{n}^2
     + \frac{\varepsilon^2}{2 L} \| \nabla \rho^n (t_n) \|^2 -
    \frac{\varepsilon^2}{2 L} \| \nabla \rho^n (t_{n - 1}) \|^2 \nn \\
    & \le D_2
    \tau \| \rho^n (t_{n - 1}) \| ^2 + D_3
    \tau \| \rho_t^{n-1} \|_{n-1}^2 + \hat{E}_p^n(\tau).
    \end{align}
   By a discrete Gr\"onwall's inequality we obtain the desired estimate
   \eqref{eq:ImexESSTerr}.
\end{proof}

\section{The semi-implicit iterative solver for the implicit scheme and superconvergence}

Now we consider the semi-implicit ESET scheme as an iterative solver for the
implicit ESET scheme: Find $h^{n, k} \in V_{M, N}^{n, u^{n-1}},$ s.t.
\begin{equation}
\int_{t_{n-1}}^{t_{n}} ( h^{n, k}_t, \phi_t ) + \varepsilon
( \nabla h^{n, k}, \nabla \phi_t ) + \frac{1}{\varepsilon}
( f (h^{n, k - 1}), \phi_t ) dt = 0, \quad \forall \; \phi
\in V^{n, 0}_{M, N }, \label{eq:itersch}
\end{equation}
for $k = 1, \ldots, N$, where $h^{n, 0} (t) = \hat{h}^n (t) : = h^{n - 1}
(t)$.

\subsection{Convergence of the iteration}

By taking difference of \eqref{eq:itersch} with $k$ and $k - 1$, and denoting
$d^{n, k}  = h^{n, k} - h^{n, k - 1}$, we have
\[ \int_{t_{n-1}}^{t_{n}} ( d_t^{n, k}, \phi_t ) + \varepsilon
( \nabla d^{n, k}, \nabla \phi_t ) + \frac{1}{\varepsilon}
( f (h^{n, k - 1}) - f (h^{n, k - 2}), \phi_t )  dt = 0,
\quad \forall \; \phi \in V^{n, 0}_{M, N}, \]
which leads to
\[ \int_{t_{n-1}}^{t_{n}} ( d^{n, k}_t, \phi_t ) + \varepsilon
( \nabla d^{n, k}, \nabla \phi_t ) + \frac{1}{\varepsilon}
( f' (\zeta^n) d^{n, k - 1}, \phi_t )  dt = 0, \quad \forall
\; \phi \in V^{n, 0}_{M, N} . \]
By taking $\phi_t = d^{n, k}_t$, we obtain
\begin{equation}
\| d_t^{n, k} \|_n^2 + \frac{\varepsilon}{2} \| \nabla d^{n, k} (t_n) \| ^2
- \frac{\varepsilon}{2} \| \nabla d^{n, k} (t_{n - 1}) \| ^2
\le
 \frac{1}{2} \| d^{n, k}_t \|^2_n
+
\frac{L^2}{2\veps^2} \| d^{n, k - 1} \|_n^2.
\end{equation}
For term
$\frac{L^2}{2 \varepsilon^2} \| d^{n, k - 1} \|_n^2$, by Lemma \ref{lem:1} and
the fact $d^{n, k-1} (t_{n - 1}) = 0$, we obtain
\[ \frac{L^2}{2 \varepsilon^2} \| d^{n, k - 1} \|_n^2 \le \frac{L^2
    \tau^2}{2 \varepsilon^2} \| d^{n, k - 1}_t \|_n^2 . \]
Combining the above two equations, we obtain
\begin{equation}
\| d^{n, k}_t \|^2_n
+ \veps\| \nabla d^{n, k}(t_n) \|^2 \le \frac{L^2
	\tau^2}{\varepsilon^2} \| d_t^{n, k - 1} \|_n^2.
\end{equation}
Thus, when $\frac{L^2 \tau^2}{\varepsilon^2} < 1$, the iteration leads to a
contraction of $\| d^{n, k}_t \|_n^2$ with rate $\frac{L^2
\tau^2}{\varepsilon^2}$. By the theorem on contracting maps, the iteration
convergences.
\begin{remark}
    The above analysis shows the convergence of the iteration under the time
 step condition $\tau\le \varepsilon/L$. Note that this is a sufficient condition,
 but may not be a necessary condition. One may develop adaptive time stepping
 method based on a posteriori error estimate of the numerical solution to
 use larger step sizes.
\end{remark}

\subsection{Uniqueness of solution to the implicit ESET scheme}

The convergence of the iterative solver in last subsection ensures the
existence of solutions to the implicit ESET scheme. Now we prove the
uniqueness.
\begin{proof}[Proof of the uniqueness part of Theorem \ref{thm:Theorem1}]
	Suppose equation \eqref{eq:SEGal} has two solutions $h^n_1$ and $h^n_2$.
	Substituting the two solutions into the equation, then taking the
	difference, and denoting by $w^n=h^n_1 - h^n_2$, we obtain
	\begin{equation*}
	    \int_{t_{n - 1}}^{t_n} ( w^n_t, v_t )
        + \varepsilon (\nabla w^n, \nabla v_t )
        + \frac{1}{\varepsilon} (f'(\zeta^n) (w^n),	v_t)  dt
        = 0, \quad \forall\, v \in V_{M,N}^{n, u^{n-1}}.
	\end{equation*}
	Then, by taking $v=w^n$, we get
	\begin{equation}
	    \| w^n_t \|^2_n	+ \frac{\varepsilon}{2} \| \nabla w^n(t_n) \|
	    - \frac{\veps}{2} \| \nabla w^n_t(t_{n-1}) \|
	    \le
	    \frac{L^2}{2\veps^2} \| w^n \|_n^2 + \frac12 \| w^n_t \|^2_n.
	\end{equation}
	Then, noticing that $w^n(t_{n-1})=0$, by using Lemma \ref{lem:1}, one
	obtain
	\begin{equation*}
	    \frac{L^2}{2\veps^2} \| w^n \|_n^2
        \le \frac{\tau^2L^2}{2\veps^2} \| w^n_t \|_n^2.
	\end{equation*}
	Combining the above results, we get
	\begin{equation}
	\left(1-\frac{\tau^2L^2}{\veps^2}\right) \| w^n_t \|^2_n
	+ {\varepsilon} \| \nabla w^n(t_n) \|
	\le 0.
	\end{equation}
	So, when $\tau < \veps/L$, we have $\| w_t^n \|_n=0$. Then, the fact
	$w^n(t_{n-1})=0$ leads to $w^n \equiv 0$.
\end{proof}

\subsection{Superconvergence of the implicit ESET scheme}

The implicit ESET scheme indeed has a better convergence rate than the
semi-implicit scheme due to superconvergence, which will be numerically
demonstrated in next section.  A lot of numerical schemes has been proved to
have superconvergence for different applications~\cite{dupont_unified_1976,
arnold_superconvergence_1979, tang_superconvergence_1992, zhou_optimal_1994,
karakashian_spacetime_1999, zhang_superconvergence_2005,
huang_superconvergence_2011, gu_superconvergence_2022,
zhang_superconvergent_2023}. Here we present a proof to the ESET scheme using a
new boosting technique, which is helpful to  better understand the
superconvergence.

To analyze the superconvergence property, we first consider the spatial
semi-discretization
\begin{equation}
(u^M_t, v) + \varepsilon ( \nabla u^M, \nabla v ) + (f(u^M), v) = 0,
\quad \forall \, v \in V_M:=P^0_M(\Omega). \label{eq:sptdis}
\end{equation}
A standard estimate will yields
\begin{equation}
    \label{eq:spterr}
    \| u - u^M \| \le C(t) \| u - \Pi_M^0 u \| \le C(t) M^{-r}.
\end{equation}
By assumption \eqref{eq:sperr-assump}, we only need to consider temporal
discretization error.
\begin{theorem}\label{thm:superconv} Let $h^n$, $u$ be the solution to
	\eqref{eq:SEGal} and \eqref{eq:SEweak}, respectively. Denote by
	$e^n(t)=h^n(t)-u(t)$. Suppose $f(u)$ satisfies
	\begin{equation} \label{eq:LipCond2}
		\max_{u} |f'(u)| < L,
		\quad
		\max_{u} |f''(u)| < L_2.
	\end{equation}
	Then we have following superconvergence result
	\begin{equation}
		\|e^n(t_n) \| \le \mathcal{O} (\tau^{2N}),
		\quad n\ge 1, \ \tau < \varepsilon/L.
	\end{equation}
\end{theorem}
\begin{proof}
1) We consider the spatial-discretized Allen--Cahn equation, with $u^M$ still
denoted by $u$:
\begin{equation} \label{eq:ACSpSemiDis}
    u _t - \varepsilon \Delta u + \frac{1}{\veps} f(u) = \chi_M(t),
\end{equation}
where $\chi_M^n$ stands for residual of spatial Galerkin projection. Let $u_r
(t)$ be a smooth solution that is close to the exact solution $u(t)$ with
error bound $\mathcal{O}(\tau^{N+1})$, e.g. $u_r |_{I_n} := \Pi^n_N \Pi^0_M u$.
We consider the linear perturbation
equation of \eqref{eq:ACSpSemiDis} near $u_r$ given below
\begin{equation}
    \label{eq:perturb}
    v_t - \varepsilon \Delta v + \frac{1}{\veps} f'(u_r) v = g (x, t).
\end{equation}
We define operator $A (t)$ as: $A (t) v := - \varepsilon \Delta v +
\frac{1}{\veps}f' (u_r(t)) v$. Then the solution of \eqref{eq:perturb} in $I_n$
can be formulated as
\begin{equation}
    v(t) = G(t_{n-1}, t) v(t_{n-1}) + \int_{t_{n-1}}^t G(s, t) g(x, s)  ds,
\label{eq:Gform}
\end{equation}
where operator $G(t_1, t_2) = e^{-\int_{t_1}^{t_2} A(s) ds}$. We first
prove the $L^2$ stability of the operator $G(s,t)$. Consider \eqref{eq:perturb}
with $g\equiv 0$, pairing the equation with $v$, we obtain
\begin{equation*}
\frac{1}{2}\frac{d}{dt}\|v \|^2 + {\veps} \|\nabla v\|^2
\le \frac{L}{\veps} \|v\|^2,
\quad t\in I_n.
\end{equation*}
Then application of Gr\"onwall's inequality leads to
\begin{equation*}
 \| v(t) \|^2 \le e^{2(t-t_{n-1})L/\veps} \| v(t_{n-1})\|^2.
\end{equation*}
Then, by \eqref{eq:Gform}, the $L^2$ operator norm of $G$ is
\begin{equation}
\| G(s,t) \| \le e^{(t-s)L/\veps}.
\end{equation}

2) Now suppose $h^n$ is the ESET solution of \eqref{eq:SEGal}, then
\begin{equation}\label{eq:Galres}
    h^n_t - \varepsilon \Delta h^n + \frac{1}{\varepsilon} f(h^n)
    = \chi^n,
\end{equation}
where $\chi^n$ stands for space-time Galerkin projection residual, which
means
\[ (\chi^n, v)_n = 0, \quad \forall \, v \in V_{M, N  - 1} . \] Taking the
difference of \eqref{eq:Galres} and \eqref{eq:ACSpSemiDis}, we have (we omit
$\chi_M(t)$ and a $(u-u_r)^2$ term, since they are smaller than
$\mathcal{O}(\tau^{2N+1})$ by assumption and standard error estimate)
\[ e^n_t - \varepsilon \Delta e^n + \frac{1}{\varepsilon} f' (u_r) e^n =
\chi^n - \frac1\veps f'' (\zeta^n) (e^n)^2
. \]
Then, by \eqref{eq:Gform}, we get
\begin{equation} \label{eq:errprop}
 e^n (t) = G (t_{n - 1}, t) e^n (t_{n - 1})
 + \int_{t_{n-1}}^{t} G
(s, t) \Bigl(\chi^n - \frac1\veps f'' (\zeta^n) (e^n(s))^2\Bigr)  ds.
\end{equation}
Taking the spatial norm on both sides, we get
\begin{align}\label{eq:mainSuper}
 \| e^n(t) \|^2
 \le 3\| G(t_{n\!-\!1}, t) e^n (t_{\!n-\!1}) \|^2
 + 3 \| \int_{t_{n\!-\!1}}^t \!\!\!G(s,t) \chi^n ds \|^2
 	+ \frac{3}{\veps^2} \| \int_{t_{n\!-\!1}}^t
     \!\!\! G(s,t) f''(\zeta^n) (e^n(s))^2 ds \|^2.
\end{align}
For the first term on the right hand side of \eqref{eq:mainSuper}, we have
\begin{equation}\label{eq:E1}
    \| G(t_{n-1}, t) e^n (t_{n-1}) \|^2
    \le e^{\frac{2(t-t_{n-1})L}{\veps}} \| e^n(t_{n-1}) \|^2.
\end{equation}
For the last term on the right hand side of \eqref{eq:mainSuper}, we have
\begin{align*}
\| \int_{t_{n-1}}^t G(s,t) f''(\zeta^n) (e^n(s))^2 ds \|^2
& \le \tau
\int_{t_{n-1}}^t  \| G(s,t) f''(\zeta^n) (e^n(s))^2 \|^2 ds \\
& \le \tau
\int_{t_{n-1}}^t  e^{\frac{2(t-s)L}{\veps}} L_2^2 \| e^n(s) \|^4_{L^4} ds \\
& \le \tau e^{\frac{2\tau L}{\veps}} L_2^2
\int_{t_{n-1}}^t K(\|\nabla e^n \| + \| e^n \| )^4 ds,\quad t\in I_n,
\end{align*}
where Sobolev embedding theorem is used in the last inequality. By Theorem
\ref{th:inpeset}, $\|\nabla e^n \| + \| e^n \| \sim \mathcal{O}(\tau^{N+1})$,
thus
\begin{align} \label{eq:E3}
\| \int_{t_{n-1}}^t G(s,t) f''(\zeta^n) (e^n(s))^2 ds \|^2
\le
L_2^2 c_K e^{\frac{2\tau L}{\veps}}  \tau^{4N+6},\quad  t\in I_n,
\end{align}
where $c_K$ is a constant independent of $\tau$.

For the second term on the right hand side of \eqref{eq:mainSuper}, we have (at
$t=t_n$)
\begin{equation}
\int_{t_{n-1}}^{t_n} G (s, t_n) \chi^n(s)  ds = \int_{t_{n-1}}^{t_n} (G
(s, t_n) - v (s)) \chi^n(s)  ds, \quad \forall \, v \in V^n_{M, N - 1} .
\label{eq:residual}
\end{equation}
Note that, since  we are working in spatial-discretized equation, $\chi^n$ can
be regarded as a vector of time, both $G$ and $v$ can be regarded as matrices.
We expand $G (s, t_n)$ at $s = t_n$. Suppose
\[ A (t) = \sum_{k = 0}^{N - 1} A^{(k)} (t_n) \frac{(t - t_n)^k}{k!} + \mathcal{O}
(\tau^N). \]
Let $Z = - \int_s^{t_{n}} A (t)  dt$, then
\[ Z = - \sum_{k = 0}^{N - 1} \frac{A^{(k)} (t_n)}{k!} \int_s^{t_{n}} (t -
    t_n)^k {dt} + \mathcal{O} (\tau^{N + 1}) = \sum_{k = 0}^{N - 1}
    \frac{A^{(k)} (t_n)}{(k + 1) !} (s - t_n)^{k + 1} + \mathcal{O}
    (\tau^{N + 1}). \] Then, $Z = \mathcal{O}(\tau)$, and
\begin{align*}
    G (s, t_n) = e^Z & = \sum_{k = 0}^{N-1} \frac{Z^k}{k!} + \mathcal{O} (\tau^{N })\\
    & = \sum_{k = 0}^{N-1} \frac{1}{k!} \left( \sum_{j = 0}^{N - 1} \frac{A^{(j)}
        (t_n)}{(j + 1) !} (s - t_n)^{j + 1} \right)^k + \mathcal{O} (\tau^{N})\\
    & = \sum_{k = 0}^{N - 1} G^{(k)} (s - t_n)^k + \mathcal{O} (\tau^N).
\end{align*}
By symmetry, $G^{(k)}$ is a symmetric diagonalizable matrix. We denote $\{
\lambda^{(k)}_i, \eta_i^{(k)} \}, i = 1, \ldots, M$ as the eigen-pairs of
$G^{(k)}$, i.e.
\[ G^{(k)} = \sum_{i = 1}^M \lambda_i^{(k)} \eta_{i}^{(k)} (\eta_i^{(k)})^T .
\] Suppose the Galerkin discretization using basis function $\psi_j (x) \in
V_M$, then the eigen-functions are given by $\varphi_i^{(k)} (x) = \sum_{j =
1}^M \eta_{ij}^{(k)} \psi_j (x)$. Let
\[ v (t, x, y) = \sum_{k = 0}^{N - 1} \sum_{i = 1}^M \lambda_i^{(k)}
\varphi^{(k)}_i (y) \varphi^{(k)}_i (x) (s - t_n)^k \quad \in V_{M, N - 1} . \]
Then
\[ G (s, t ; x, y) - v (t, x, y) \sim \mathcal{O} (\tau^{N}) . \] And it is not
hard to show that the residual error has a bound $\| \chi^n \| \sim \mathcal{O}
(\tau^N)$. So we have
\begin{equation}
\label{eq:E2}
\| \int_{t_{n-1}}^t \!\!\!G(s,t) \chi^n ds \|^2 \le \mathcal{O} (\tau^{4N+1}).
\end{equation}
Combine \eqref{eq:mainSuper}, \eqref{eq:E1}, \eqref{eq:E3} and \eqref{eq:E2},
then use a discrete Gr\"onwall's inequality, we obtain the superconvergence result.
\end{proof}

\begin{remark}\label{rmk:2} As discussed in Remark \ref{rmk:1}, the standard
    double-well potential \eqref{eq:double-wellF} doesn't satisfy condition
    \eqref{eq:LipCond}. But the truncated quadratic growth potential
    \eqref{eq:truncatedF-M}
    and \eqref{eq:truncatedF-M1} satisfies the Lipschitz condition
    \eqref{eq:LipCond} with $L =3M^2-1$ and $L=2$, respectively. Actually, the
    second condition in \eqref{eq:LipCond2} is satisfied almost everywhere with
    $L_2=6M$ and $L_2=6$, respectively for potential \eqref{eq:truncatedF-M}
    and \eqref{eq:truncatedF-M1}. To have smaller $L$ and $L_2$, we suggest to
    use \eqref{eq:truncatedF-M1} for the Allen--Cahn equation.
\end{remark}

\section{Numerical results}

In this section, we numerically verify the stability and accuracy of
proposed schemes.

To test the numerical scheme, we first solve (\ref{eq:AC}) in a one-dimensional
domain $\Omega = [-1,1]$. For the homogeneous Dirichlet boundary condition
$u\lvert_{\partial\Omega} = 0$, we take a reference solution
\begin{equation} \label{eq:uexDBC}
    u_\text{ref} = \tanh(\tfrac{x-t}{\varepsilon})
                    -\frac{1}{2}[(x+1)\tanh(\tfrac{1-t}{\varepsilon})
                        +(1-x)\tanh(\tfrac{-1-t}{\varepsilon})].
\end{equation}
Neumann boundary condition can be achieved similarly. Notice that the reference
solution does not satisfy Allen--Cahn equation (\ref{eq:AC}). We add an
external force term into (\ref{eq:AC}) to make \eqref{eq:uexDBC} an exact
solution.

\begin{figure}[htbp]
    \centering
    \includegraphics[width=0.44\textwidth]{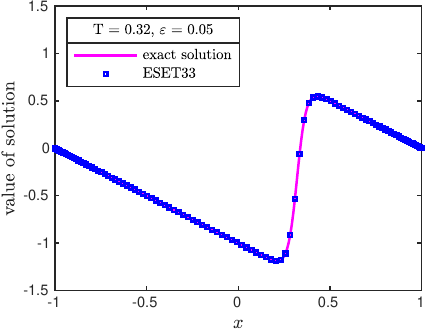}\quad
    \includegraphics[width=0.44\textwidth]{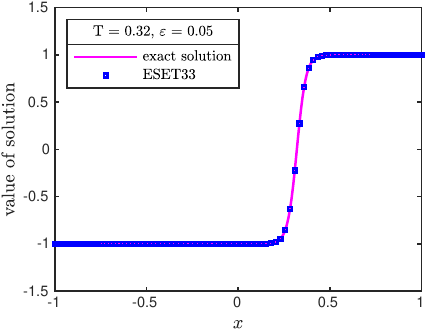}
    \caption{\label{fig:1}Testing ESET with Dirichlet(left) and Neumann(right)
    boundary condition, respectively. The scheme parameters used: $M =
    255, N=3$, $\tau = 0.01$, $S=0$.}
\end{figure}

First, we take $T = 0.32,M = 255,\varepsilon = 0.05$, $S=0$, $\tau =
0.01$ and use two different boundary conditions to test the stability and
accuracy of the proposed schemes. The numerical results are given in Figure
\ref{fig:1}, from which we see the proposed scheme give good numerical
solutions. Here ESET33 means using semi-implicit ESET scheme with $N=3$ (first
3 in ESET33) as an iterative solver for the implicit ESET scheme using $3$
(second 3 in ESET33) iterations. From the analysis given in last section,
we know that each iteration will increase the order of the numerical scheme by
1, until the maximum order (which is 6 for the implicit scheme with $N=3$) is
achieved. Using more than 3 iterations can further increase the accuracy of
numerical solutions a little bit, but the convergence order can't be improved.

\begin{figure}[htbp]
    \centering
    \includegraphics[width=0.44\textwidth]{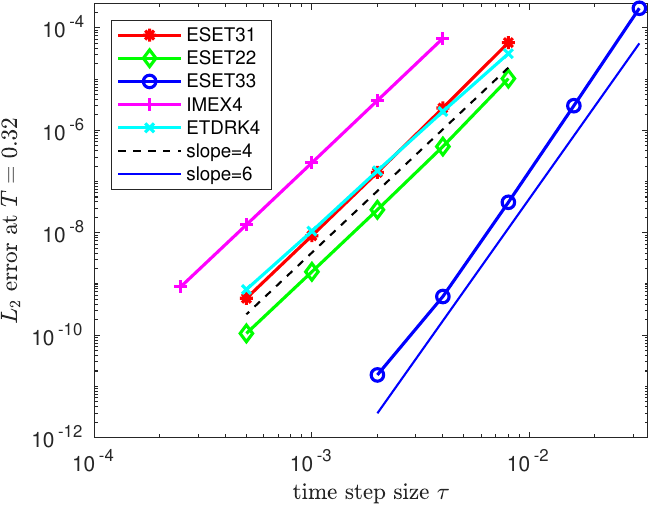}\quad
    \includegraphics[width=0.44\textwidth]{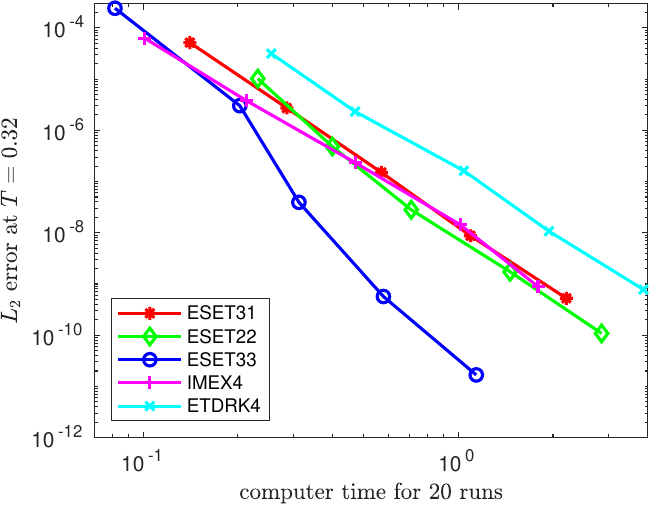}
    \caption{\label{fig:AccCost}Accuracy and efficiency of ESET, IMEX4, and
    ETDRK4 schemes. The common parameters used: $\veps=0.05$, $M = 350$,
    $S=0$.}
\end{figure}

We then test convergence rates and make comparison with well-known schemes.
Several commonly used fourth order time marching schemes are investigated by
Kassam and Trefethen~\cite{kassam_fourthorder_2005}, they found that the fourth
order implicit-explicit backward differentiation formula (IMEX4), and the
fourth order exponential time-differencing Runge-Kutta method (ETDRK4) give the
best performance results. We compare our scheme with these two schemes. The
convergence rates of the proposed implicit (ESET22, ESET33) and semi-implicit
(ESET31) schemes together with IMEX4 and ETRK4 schemes are presented in the
left plot of Figure \ref{fig:AccCost}. The results of numerical accuracy versus
computer time for these schemes are given in the right plot of Figure
\ref{fig:AccCost}. We find that ESET31, ESET22, IMEX4, ETDRK4 all have fourth
order convergence, but ESET22 has the smallest error constant and is the most
efficient one among these fourth order schemes in terms of computational cost,
due to superconvergence. When high accuracy is needed, the 6th order ESET33
scheme out-perform all the fourth order schemes significantly.

\begin{figure}[htbp]
    \centering
    \includegraphics[width=0.44\textwidth]{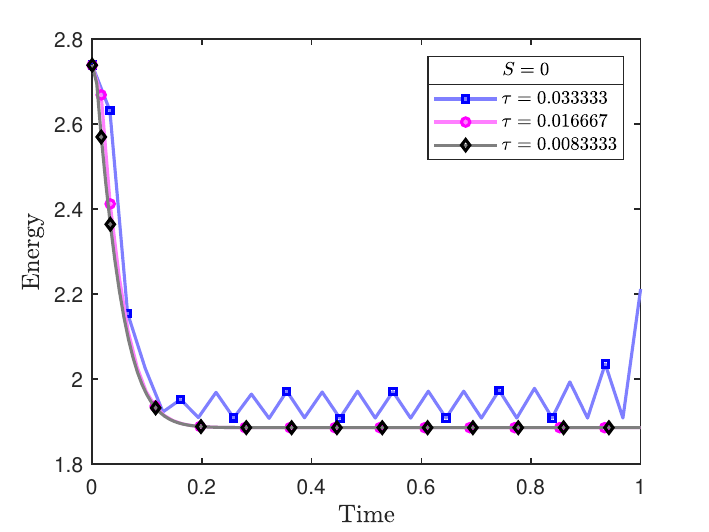}\quad
    \includegraphics[width=0.44\textwidth]{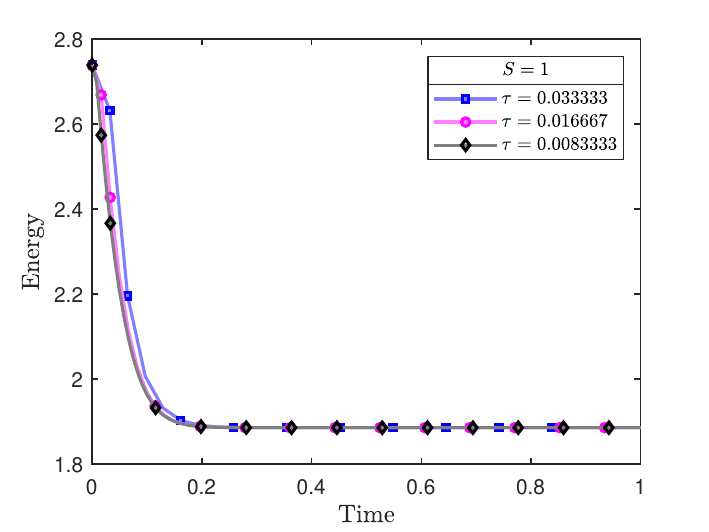}
    \caption{\label{fig:elaw}The energy dissipation of ESET31 scheme
    \eqref{eq:SEGalLinear} with different stabilization constant $S$ for
    $\varepsilon=0.08$ and $M=350$. }
\end{figure}

Figure \ref{fig:elaw} presents the energy dissipation of ESET31 scheme
\eqref{eq:SEGalLinear} with different time step sizes and stabilization
constants. We find that shortening the time step can stabilize the energy
curve, and the stability constant is also helpful.

\begin{figure}[htbp]
    \centering
    \includegraphics[width=0.44\textwidth]{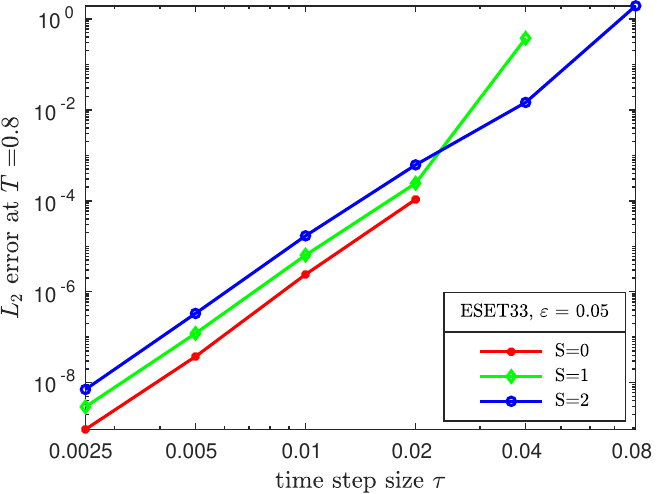}\quad
    \includegraphics[width=0.44\textwidth]{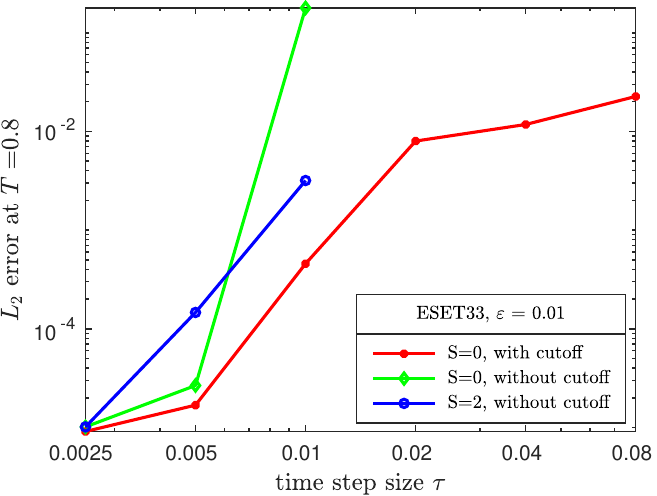}
    \caption{\label{fig:StabCutoff}Effects on stability and accuracy of the
    stabilization constant and cut-off operation to maintain maximum principle.
    $M=350$. In this figure, the lines break at a time step size when the
    corresponding numerical scheme encounters blow-ups.}
\end{figure}

In Figure \ref{fig:StabCutoff}, we investigate the effect on accuracy of
stabilization and a simple cut-off operation to maintain maximum
principle~\cite{li_arbitrarily_2020, yang_arbitrarily_2022}. From this figure,
we see that the stabilization constant $S$ can improve the stability when time
step size is large, but will increase the numerical error a little bit when
time step size is small and the stability is not an issue. This suggests that
one should adjust the stability constant according to time step
size~\cite{wang_energystable_2019}. On the other side, using cut-off increases
both stability and accuracy, see the right plot of Figure \ref{fig:StabCutoff}.

\begin{figure}[htbp]
    \centering
    \includegraphics[width=0.85\textwidth]{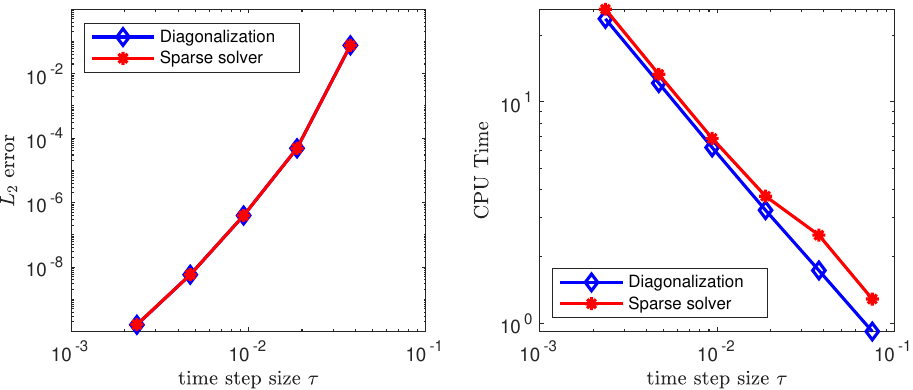}
    \caption{\label{fig:solver_method}Comparison of accuracy and efficiency
    between the diagonalization method and direct sparse solver
    of the ESET42 scheme with $M=350, \veps=0.08, T=1.2$.}
\end{figure}
In Figure \ref{fig:solver_method}, we compare the performance of the direct
sparse solver and the diagonalization approach in a ESET42 scheme. We see two
solution methods give the same numerical solutions, and the diagonalization
method use slightly less computer time.

\medskip
Next we solve the conservative Allen--Cahn
equation~\cite{rubinstein_nonlocal_1992} in a two-dimensional domain $\Omega =
[-1,1]^2$:
\begin{equation}\label{eq:CAC}
    \frac{\partial u}{\partial t}(x,t) - \varepsilon\Delta u(x,t)
    = - \tfrac{1}{\varepsilon}f(u(x,t)) + \alpha(t),   \quad x\in \Omega,
\end{equation}
where
\begin{equation}
\alpha(t) = \frac{1}{|\Omega|}{\int_{\Omega} f(u(x,t))dx},
\end{equation}
is a non-local term to conserve total mass
\begin{equation*}
\int_{\Omega}u(x,t)dx = \int_\Omega u(x,0) dt.
\end{equation*}
The numerical scheme for \eqref{eq:CAC} is almost the same as for the standard
Allen--Cahn equation, since $\alpha$ is a function depends only on $f(u)$.

More precisely, the semi-implicit ESET scheme for \eqref{eq:CAC} reads:
Find $h^n \in V_{M, N}^{n, u^{n - 1}},$ s.t.
\begin{equation}\label{eq:SEGalCAC}
    \int_{t_{n - 1}}^{t_n} (h^n_t, v_t)
    + \varepsilon (\nabla h^n, \nabla v_t )
    + (\frac{1}{\varepsilon} f(\hat{h}^n) - \hat{\alpha}(t), v_t)  dt
    = 0, \quad \forall\, v \in V_{M,N}^{n, u^{n-1}},
\end{equation}
where $\hat{\alpha}(t) = \frac{1}{|\Omega|}{\int_{\Omega}
f(\hat{h}^n(x,t))dx}$. By taking test function $v_t\equiv 1$ in
\eqref{eq:SEGalCAC}, we obtain the discrete mass conservation automatically.

\begin{figure}[htbp]
    \centering
    \includegraphics[width=0.22\textwidth]{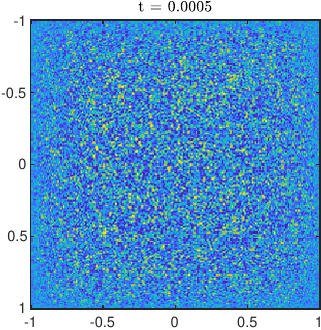}
    \includegraphics[width=0.22\textwidth]{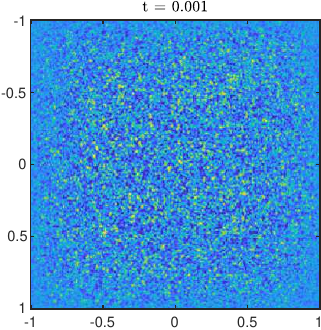}
    \includegraphics[width=0.22\textwidth]{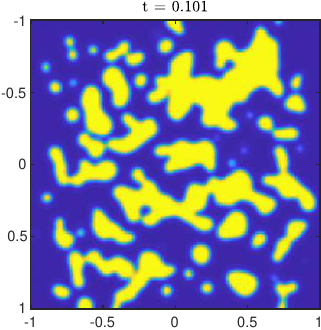}
    \includegraphics[width=0.22\textwidth]{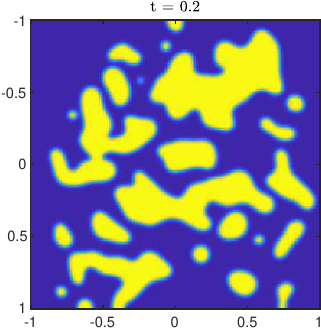}\\
        \centering
    \includegraphics[width=0.22\textwidth]{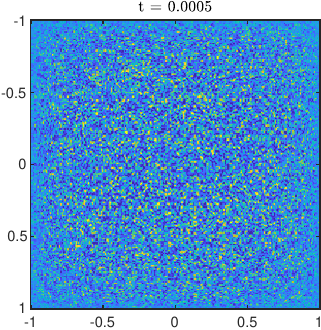}
    \includegraphics[width=0.22\textwidth]{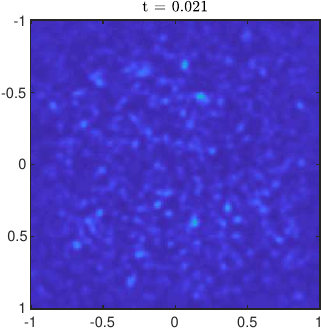}
    \includegraphics[width=0.22\textwidth]{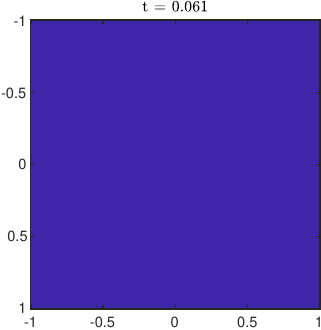}
    \includegraphics[width=0.22\textwidth]{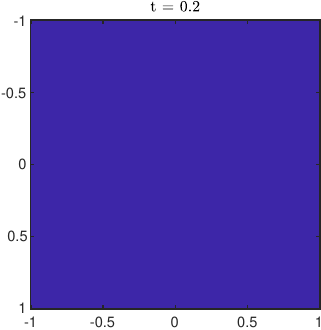}
    \caption{\label{fig:ACrandomIV}Solution snapshots of conservative
    Allen–Cahn equation ({\em top}) and standard Allen--Cahn equation ({\em
    bottom}) for random initial condition \eqref{eq:random_IC} by the
    ESET22 scheme. $\veps=0.01$, $M=280^2$, $T=0.2$, $S=0$.}
\end{figure}

\begin{figure}[htbp]
    \centering
    \includegraphics[width=0.44\textwidth]{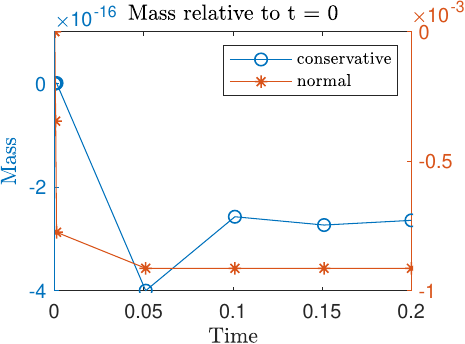}\quad
    \includegraphics[width=0.44\textwidth]{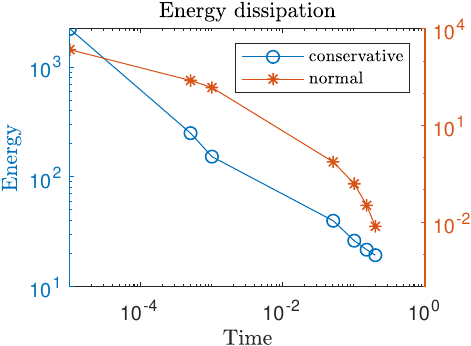}
    \caption{\label{fig:massEnergy}Mass changes (left) and energy dissipation
    (right) of the conservative Allen--Cahn equation and standard (normal) Allen--Cahn
    equation for random initial condition \eqref{eq:random_IC} by the
    ESET22 scheme. $\veps=0.01$, $M=280^2$, $T=0.2$, $S=0$.}
\end{figure}

We first use ESET22 to simulate a case with $\varepsilon = 0.01$ and random
initial value. The initial solution takes uniform random number $U(0,1)$ at
each tensor product Legendre-Gauss quadrature point $(x_i, y_j)$:
\begin{equation}\label{eq:random_IC}
    h^0(x_i,y_j) \sim U(0,1).
\end{equation}
We use $M=280^2$ spatial bases, and set initial time step size to be $\tau =
10^{-5}$, then increase the time step to $10^{-3}$ after calculating 99 steps.
Figure \ref{fig:ACrandomIV} presents the snapshots of solutions of conservative
Allen--Cahn equation together with standard Allen--Cahn equation with the same
initial condition for comparison. We observe that for conserved case, phase
separation happened in a very short time, while one phase gets dominant in the
non-conserved case. The corresponding mass conservation and energy dissipation
are presented in Figure \ref{fig:massEnergy}, from which we observe that the
total mass is kept up to machine accuracy for the conservative Allen--Cahn
equation, but no mass conservation for the standard Allen--Cahn equation,
which are consistent with the results shown in Figure
\ref{fig:ACrandomIV}. We observe that both cases dissipate energy.
\pagebreak

\begin{figure}[htbp]
    \centering
    \includegraphics[width=0.22\textwidth]{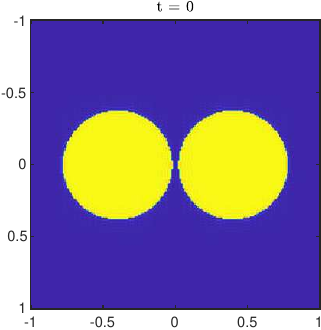}
    \includegraphics[width=0.22\textwidth]{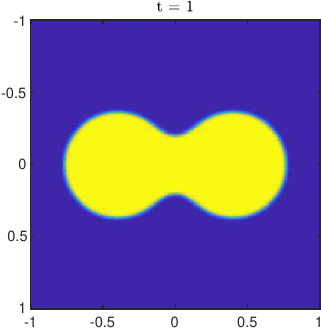}
    \includegraphics[width=0.22\textwidth]{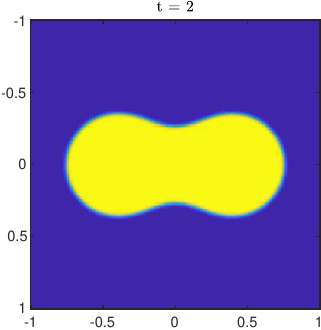}
    \includegraphics[width=0.22\textwidth]{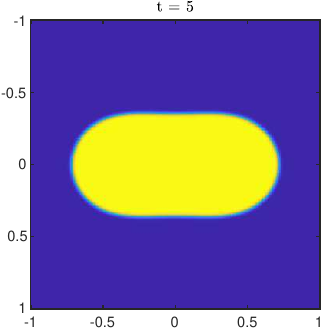}
    \includegraphics[width=0.22\textwidth]{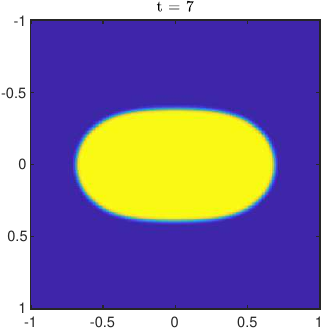}
    \includegraphics[width=0.22\textwidth]{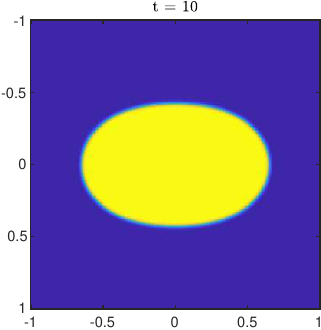}
    \includegraphics[width=0.22\textwidth]{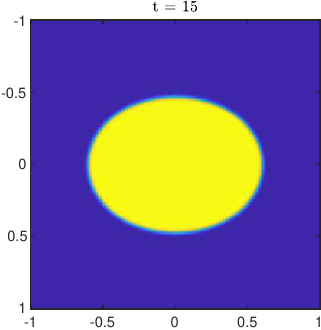}
    \includegraphics[width=0.22\textwidth]{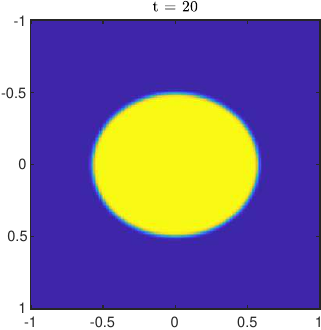}
    \caption{\label{fig:BallsSnapCAC}Snapshots of conservative Allen--Cahn
    equation for the drop coalescence test using ESET22 scheme with
    $\tau=5\times 10^{-3}$, $M=280^2$, $\veps=0.01$, $T=20$.}
\end{figure}

\begin{figure}[htbp]
    \centering
    \includegraphics[width=0.22\textwidth]{bubble_init.pdf}
    \includegraphics[width=0.22\textwidth]{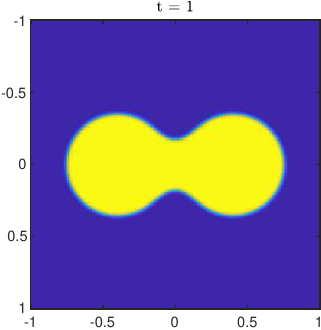}
    \includegraphics[width=0.22\textwidth]{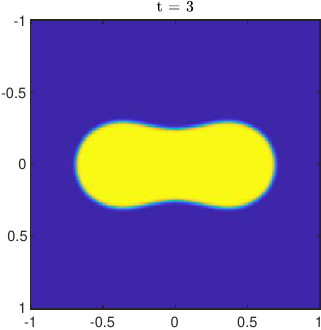}
    \includegraphics[width=0.22\textwidth]{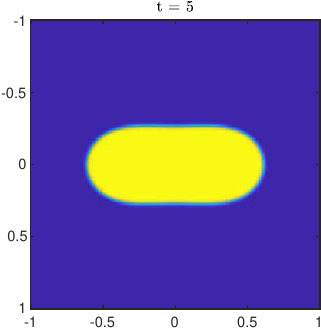}
    \includegraphics[width=0.22\textwidth]{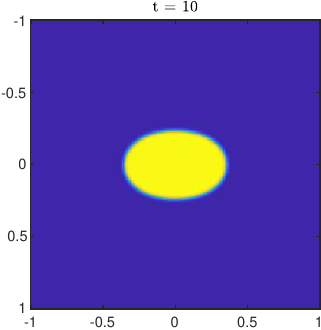}
    \includegraphics[width=0.22\textwidth]{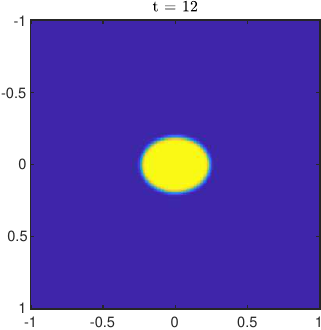}
    \includegraphics[width=0.22\textwidth]{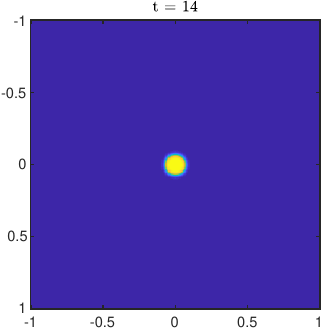}
    \includegraphics[width=0.22\textwidth]{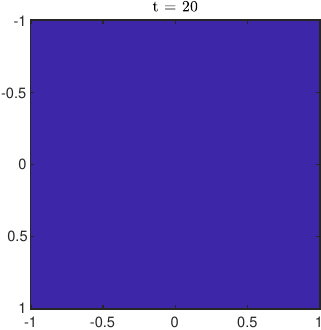}
    \caption{\label{fig:BallsSnapSAC}Snapshots of standard Allen–Cahn equation
    for the drop coalescence test using ESET22 scheme with $\tau=5\times
    10^{-3}$, $M=280^2$, $\veps=0.01$, $T=20$.}
\end{figure}

\begin{figure}[htbp]
    \centering
    \includegraphics[width=0.42\textwidth]{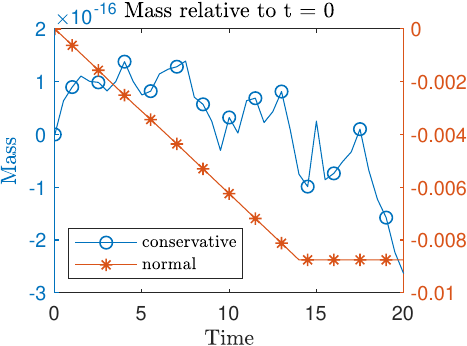}\quad
    \includegraphics[width=0.42\textwidth]{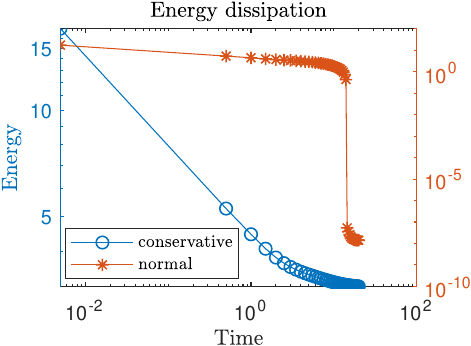}
    \caption{\label{fig:BallMassEnergy}Mass changes (left) and energy
    dissipation (right) of conservative and standard (normal) Allen--Cahn equations
    for the drop coalescence test obtained by ESET22 scheme with $\tau=5\times
    10^{-3}$, $M=280^2$, $\veps=0.01$, $T=20$.}
\end{figure}

Next, we simulate drop coalescence with initial state contains two balls
centered at $(0.4, 0)$ and $(-0.4,0)$ with radius 0.38. We use same solver
ESET22 to solve the standard and conservative Allen--Cahn equations. The
results are presented in Figure \ref{fig:BallsSnapCAC}, \ref{fig:BallsSnapSAC}
and \ref{fig:BallMassEnergy}. We observe that in the conserved case,  two balls
merge into one larger ball, the total mass is kept; but in standard case, the
two balls first merge then disappear at an almost constant rate of total
mass diminishing. Both cases dissipate energy.

\section{Concluding remarks}
Based on an energetic variational formulation, we proposed two efficient
spectral-element time marching methods for nonlinear gradient systems: one
implicit ESET scheme and one semi-implicit ESET scheme. The semi-implicit ESET
scheme leads to linear systems with constant coefficients, thus can be
efficiently solved. The implicit ESET scheme has superconvergence property, it
can be efficiently solved by using the semi-implicit scheme as an iterative
solver. There are several advantages of the proposed schemes: 1) They keep mass
conservation accurately if the continuous equation conserves mass. 2) They keep
the energy dissipation with a time step size related to physical time scale. 3)
They are high order accurate and efficient. In particular, when using the
semi-implicit scheme as an iterative solver for the fully implicit scheme,
superconvergence can be achieved, resulting a more efficient method than
existing ones. 4) The diagonalization solution procedure allows the method to
be used for large scale problems with parallel computing in time.

In this paper, we only considered the Allen--Cahn equations. But the proposed
methods can be used for other parabolic type nonlinear systems. However,
application to $H^{-1}$ gradient systems, e.g. the Cahn--Hilliard equation is
not trivial, direct extension without a good stabilization leads to a
scheme with time step size restriction of order $\mathcal{O}(\veps^3)$. How to
design unconditionally stable linear ESET schemes for Allen--Cahn type and
Cahn--Hilliard type equations deserves a further study.

\section*{CRediT authorship contribution statement}
Shiqin Liu: Conceptualization, Formal analysis, Investigation, Methodology,
Writing. Haijun Yu: Conceptualization, Formal analysis, Investigation,
Methodology, Writing, Funding acquisition.

\section*{Declaration of competing interest}
We declare that we have no financial and personal relationships with other
people or organizations that can inappropriately influence our work. There is
no professional or other personal interest of any nature or kind in any
product, service and/or company that could be construed as influencing the
position presented in, or the review of, the manuscript entitled.

\section*{Acknowledgments}

This work was partially supported by the National Natural Science Foundation of
China (grant number 12171467, 12161141017). The computations were partially
done on the high-performance computers of the State Key Laboratory of
Scientific and Engineering Computing, Chinese Academy of Sciences.

\bibliography{eset_paper}
\bibliographystyle{elsarticle-num}

\end{document}